\documentclass[reqno,10pt,centertags,draft]{amsart}
\usepackage{amsmath,amsthm,amscd,amssymb,eufrak,latexsym,upref}

\makeatletter
\def\theequation{\@arabic\c@equation}

\newcommand{\diag}{\operatorname{diag}}

\newcommand{\bbN}{{\mathbb{N}}}
\newcommand{\bbR}{{\mathbb{R}}}

\newcommand{\bbC}{{\mathbb{C}}}

\newcommand{\cA}{{\mathcal A}}
\newcommand{\cB}{{\mathcal B}}
\newcommand{\cC}{{\mathcal C}}

\newcommand{\cF}{{\mathcal F}}
\newcommand{\cG}{{\mathcal G}}
\newcommand{\cH}{{\mathcal H}}
\newcommand{\cI}{{\mathcal I}}

\newcommand{\cL}{{\mathcal L}}
\newcommand{\cM}{{\mathcal M}}
\newcommand{\cN}{{\mathcal N}}

\newcommand{\cP}{{\mathcal P}}
\newcommand{\cQ}{{\mathcal Q}}

\newcommand{\cU}{{\mathcal U}}

\newcommand{\gC}{\mathfrak{C}}
\newcommand{\gh}{\mathfrak{h}}
\renewcommand{\gg}{\mathfrak{g}}

\newcommand{\gR}{\mathfrak{R}}

\newcommand{\no}{\nonumber}
\newcommand{\lb}{\label}
\newcommand{\f}{\frac}
\newcommand{\ul}{\underline}
\newcommand{\ol}{\overline}

\newcommand{\wti}{\widetilde  }

\newcommand{\oh}{o}

\newcommand{\loc}{\text{\rm{loc}}}
\newcommand{\spec}{\text{\rm{spec}}}
\newcommand{\rank}{\text{\rm{rank}}}

\newcommand{\dom}{\text{\rm{dom}}}
\newcommand{\ess}{\text{\rm{ess}}}
\newcommand{\p}{\text{\rm{p}}}
\newcommand{\ac}{\text{\rm{ac}}}
\newcommand{\singc}{\text{\rm{sc}}}

\newcommand{\AC}{\text{\rm{AC}}}
\newcommand{\bi}{\bibitem}
\newcommand{\hatt}{\widehat}

\newcommand{\skdv}{\text{\rm{s-KdV}}}

\renewcommand{\Re}{\text{\rm Re}}
\renewcommand{\Im}{\text{\rm Im}}
\renewcommand{\ln}{\text{\rm ln}}

\numberwithin{equation}{section}

\newtheorem{theorem}{Theorem}[section]
\newtheorem{lemma}[theorem]{Lemma}
\newtheorem{corollary}[theorem]{Corollary}
\newtheorem{hypothesis}[theorem]{Hypothesis}
\newtheorem{definition}[theorem]{Definition}
\newtheorem{remark}[theorem]{Remark}

\newcommand{\abs}[1]{\lvert#1\rvert}

\begin{document}

\title[Extensions of Theorems of Borg and Hochstadt]{Matrix-Valued
Generalizations  of the Theorems of  Borg and Hochstadt}

\author[E.\ D.\ Belokolos]{Eugene\ D.\ Belokolos}
\address{Institute of Magnetism, National Academy of Sciences,
Ukraine, Vernadsky Str.\ 36B, Kiev--142, 252142, Ukraine}
\email{bel@imag.kiev.ua}

\author[F.\ Gesztesy]{Fritz\ Gesztesy}
\address{Department of Mathematics,
University of Missouri, Columbia, MO 65211, USA}
\email{fritz@math.missouri.edu}
\urladdr{http://www.math.missouri.edu/people/fgesztesy.html}

\author[K.\ A.\ Makarov]{Konstantin\ A.\ Makarov}
\address{Department of Mathematics,
University of Missouri, Columbia, MO 65211, USA}
\email{makarov@math.missouri.edu\newline
\indent{\it URL:}
http://www.math.missouri.edu/people/kmakarov.html}

\author[L.\ A.\ Sakhnovich]{\mbox{\,Lev\ A.\ Sakhnovich}}
\address{735 Crawford Ave., Brooklyn, NY 11223, USA}
\email{Lev.Sakhnovich@verizon.net}

\dedicatory{Dedicated with great pleasure to Jerry 
Goldstein and Rainer Nagel \\ on the occasion of their 60th birthdays}

\date{December, 2001}
\thanks{Research of the first and second author was supported in part
by the CRDF grant UM1-325.}

\begin{abstract}
We prove a generalization of the well-known theorems by Borg and
Hochstadt for periodic self-adjoint Schr\"odinger operators
without a spectral gap, respectively, one gap in their spectrum, in the
matrix-valued context. Our extension of the theorems of Borg and
Hochstadt replaces the periodicity condition of the potential by the more
general property of being reflectionless (the resulting potentials then
automatically turn out to be periodic and we recover Despr\'es' matrix
version of Borg's result). In addition, we assume the spectra to have
uniform maximum multiplicity (a condition automatically fulfilled in the
scalar context considered by Borg and Hochstadt). Moreover, the
connection with the stationary matrix KdV hierarchy is established. 

The methods employed in this paper rely on matrix-valued Herglotz
functions, Weyl--Titchmarsh theory, pencils of matrices, and basic inverse
spectral theory for matrix-valued Schr\"odinger operators.
\end{abstract}

\maketitle

\section{Introduction} \lb{s1}

In a previous paper, \cite{GS02}, two of us constructed a class of
self-adjoint $m\times m$ matrix-valued Schr\"odinger operators 
$H(\Sigma_n)=-d^2/dx^2 \cI_m+\cQ(\Sigma_n,\cdot)$ in $L^2(\bbR)^{m\times
m}$,
$m\in\bbN$, with prescribed absolutely continuous finite-band spectrum
$\Sigma_n$ of the type
\begin{equation}
\Sigma_n=\Bigg\{\bigcup_{j=0}^{n-1} [E_{2j},E_{2j+1}]\Bigg\}\cup
[E_{2n},\infty), \quad n\in\bbN_0  \lb{1.1}
\end{equation}
of uniform spectral multiplicity $2m$. Here 
\begin{equation}
\{E_\ell\}_{0\leq \ell\leq 2n}\subseteq \bbR, \; n\in\bbN, \text{ with
$E_\ell<E_{\ell+1}$, $0\leq \ell\leq 2n-1$,} \lb{1.2}
\end{equation}
and hence $H(\Sigma_n)$ satisfies
\begin{equation}
\spec(H(\Sigma_n))=\Sigma_n. \lb{1.3}
\end{equation}

Throughout this paper all matrices will be considered over
the field of complex numbers $\bbC$, and the corresponding linear space of
$k\times\ell$ matrices will be denoted by $\bbC^{k\times\ell}$,
$k,\ell\in\bbN$.  Moreover, $\cI_k$ denotes the identity matrix in
$\bbC^{k\times k}$ for $k\in\bbN$, $\cM^*$ the adjoint (i.e., complex
conjugate transpose), $\cM^t$ the transpose of a matrix $\cM$,
$\diag(m_1,\dots,m_k)\in\bbC^{k\times k}$ a diagonal $k\times k$
matrix, and $\AC_{\loc}(\bbR)$ denotes the set of locally absolutely
continuous functions on $\bbR$. The spectrum, point spectrum (the
set of eigenvalues), essential spectrum, absolutely continuous spectrum,
and singularly continuous spectrum of a self-adjoint linear operator
$T$ in a separable complex Hilbert space are denoted by $\spec(T)$,
$\spec_{\p}(T)$, $\spec_{\ess}(T)$, $\spec_{\ac}(T)$, and
$\spec_{\singc}(T)$, respectively.

The constructed matrix potentials $\cQ_{\Sigma_n}$ (resp.,
$H(\Sigma_n)$) turns out to be {\it reflectionless} in the sense 
discussed in \cite{CGHL00}, \cite{GT00}, and \cite{KS88}
(cf.\ \cite{Cr89}, \cite{DS83}, \cite{KK88}, \cite{SY95} in the scalar
context $m=1$), that is, the half-line Weyl--Titchmarsh matrices
$\cM_{\pm}(\Sigma_n,z,x)$ associated with $H(\Sigma_n)$, the half-lines
$[x,\pm\infty)$, and a Dirichlet boundary condition at $x\in\bbR$,
satisfy 
\begin{align}
&\lim_{\varepsilon\downarrow 0}\cM_{+}(\Sigma_n,\lambda+i\varepsilon,x)
=\lim_{\varepsilon \downarrow
0}\cM_{-}(\Sigma_n,\lambda-i\varepsilon,x), \lb{1.4} \\
& \hspace*{1.3cm} \lambda\in \bigcup_{j=0}^{n-1} (E_{2j},E_{2j+1}) \cup
(E_{2n},\infty), \; x\in\bbR. \no
\end{align}
Especially, $\cM_{+}(\Sigma_n,\cdot,x)$ is the analytic continuation
of $\cM_{-}(\Sigma_n,\cdot,x)$ through the set $\Sigma_n$, and vice
versa. In other words, $\cM_{+}(\Sigma_n,\cdot,x)$ and
$\cM_{-}(\Sigma_n,\cdot,x)$ are the two branches of an analytic
matrix-valued function $\cM(\Sigma_n,\cdot,x)$ on the two-sheeted
Riemann surface of $\big(\prod_{\ell=0}^{2n} (z-E_\ell)\big)^{1/2}$. The
reflectionless property \eqref{1.4} then implies the absolute
continuity of the spectrum of $H_{\Sigma_n}$ and its uniform (maximal)
multiplicity $2m$.

In this sequel of paper \cite{GS02}, we focus on the two special cases
$n=0$ and $n=1$ and derive matrix-valued extensions of the
well-known theorems by Borg \cite{Bo46} and Hochstadt \cite{Ho65},
respectively. Before  describing our principal new results and the
contents of each section, we briefly recall the classical results by Borg
and Hochstadt in the scalar case $m=1$. 

In 1946 Borg \cite{Bo46} proved, among a variety of other inverse
spectral theorems, the following result.

\begin{theorem}[\cite{Bo46}] \lb{t1.1}
Let $\Sigma_0=[E_0,\infty)$ for some $E_0\in\bbR$ and
$q_0\in L^1_{\loc} (\bbR)$ be real-valued and periodic.
Suppose that $h_0=-\f{d^2}{dx^2}+q_0$ is the associated
self-adjoint  Schr\"odinger operator in $L^2(\bbR)$ $($cf.\ \eqref{1.9}
for $m=1$$)$ and assume that
\begin{equation}
\spec (h_0)=\Sigma_0. \lb{1.5}
\end{equation}
Then
\begin{equation}
q_0(x)=E_0 \text{ for a.e. $x\in\bbR$}. \lb{1.6}
\end{equation}
\end{theorem}
Traditionally, uniqueness results such as Theorem \ref{t1.1} are called
Borg-type theorems. (However, this terminology is not uniquely adopted 
and hence a bit unfortunate. Indeed, inverse spectral results on finite
intervals recovering the  potential coefficient(s) from several spectra,
were also pioneered by Borg in his celebrated paper \cite{Bo46}, and hence
are also coined Borg-type theorems in the literature, see, e.g., 
\cite[Sect.\ 6]{Ma94}.) Actually, Borg assumed $q_0\in
L^2_{\loc}(\bbR)$, but that seems a minor detail. 

\begin{remark} \lb{r1.2} ${}$ \\ 
$(i)$ A closer examination of the short proof of $($an extension of\,$)$ 
Theorem \ref{t1.1} provided in \cite{CGHL00} shows that periodicity  of
$q_0$ is not the point for the uniqueness result \eqref{1.6}. 
The key ingredient $($besides $\spec (h_0)=[E_0,\infty)$ and
$q_0$  real-valued\,$)$ is clearly the fact that $q_0$
is reflectionless in the sense of \eqref{1.4}. \\
$(ii)$ Real-valued periodic potentials are known to satisfy \eqref{1.4}, 
but so are certain classes of real-valued quasi-periodic and
almost-periodic potentials $q_0$ $($see, e.g., \cite{Cr89},
\cite{DS83}, \cite{Ko84}, \cite{Ko87a}, \cite{Ko87b}, \cite{KK88},
\cite{KS88}, \cite{SY95}$)$. In particular, the class of real-valued
algebro-geometric  finite-gap potentials $q_0$ $($a subclass of
the set of real-valued quasi-periodic  potentials$)$ is a prime example
satisfying \eqref{1.4} without necessarily being periodic. \\
$(iii)$ We note that real-valuedness of $q_0$ is an essential 
assumption in Theorem~\ref{t1.1}. Indeed, it is well-known that
$q(x)=\exp(ix)$, $x\in\bbR$, leads to the half-line spectrum $[0,\infty)$.
A detailed treatment of a class of examples  of this type can be found in
\cite{Ga80}, \cite{Ga80a}, \cite{GU83}, \cite{PT88}, \cite{PT91}.
Moreover, the example  of complete exponential localization of the
spectrum of a  discrete Schr\"odinger operator with a quasi-periodic 
real-valued potential having two basic frequencies and no gaps in its
spectrum \cite{CS89} illustrates the importance of the reflectionless
property of $q_0$ in Theorem~\ref{t1.1}.
\end{remark}

Next we recall Hochstadt's theorem \cite{Ho65} from 1965.

\begin{theorem}[\cite{Ho65}] \lb{t1.3}
Let $\Sigma_1=[E_0,E_1]\cup [E_2,\infty)$ for some $E_0<E_1<E_2$ and 
$q_1\in L^1_{\loc} (\bbR)$ be real-valued and periodic.
Suppose that $h_1=-\f{d^2}{dx^2}+q_1$ is the
associated self-adjoint  Schr\"odinger operator in $L^2(\bbR)$
$($cf.\ \eqref{1.9} for $m=1$$)$ and assume that
\begin{equation}
\spec (h_1)=\Sigma_1. \lb{1.7}
\end{equation}
Then 
\begin{equation}
q_1(x)=C_0 +2\wp(x+\omega_3+\alpha) 
\text{ for some $\alpha\in\bbR$ and a.e.\ $x\in\bbR$.} \lb{1.8}
\end{equation}
\end{theorem}
Here $\wp(\cdot)=\wp(\cdot;\omega_1;\omega_3)$ denotes the elliptic
Weierstrass function with half-periods $\omega_1>0$ and $-i\omega_3>0$ 
(cf.\ \cite[Ch.\ 18]{AS72}). 

\begin{remark} \lb{r1.4} Again it will turn out that periodicity 
of $q_1$ is not the point for the uniqueness result \eqref{1.8}.
The key ingredient $($besides $\spec (h_1)=\Sigma_1$ and
$q_1$  real-valued\,$)$ is again  the fact that $q_1$ is
reflectionless in the sense of \eqref{1.4}. Similarly, Remarks
\ref{r1.2}\,$(ii),(iii)$ apply of course in the present context.
\end{remark}

The principal results of this paper then read as follows.

\begin{theorem} \lb{t1.5} Let $m\in\bbN$, suppose 
$\cQ_\ell=\cQ_\ell^* \in
L^1_{\loc}(\bbR)^{m\times m}$ and assume that the differential
expressions $-\cI_m\f{d^2}{dx^2}+\cQ_\ell$, $\ell=0,1$, are in
the limit  point case at $\pm\infty$. Define the self-adjoint
Schr\"odinger operators $H_\ell$ in $L^2(\bbR)^{m\times m}$ 
\begin{align}
&H_\ell=-\cI_m \f{d^2}{dx^2}+\cQ_\ell, \quad
\ell=0,1, \lb{1.9} \\
&\dom(H_\ell)=\{g\in L^2(\bbR)^m \mid g,g^\prime\in
\AC_{\loc}(\bbR)^m;\, (-g^{\prime\prime}+\cQ_\ell g)\in
L^2(\bbR)^m\} \no
\end{align}
and assume that $\cQ_\ell$ is reflectionless $($cf.\ \eqref{1.4}$)$. \\
$(i)$ Let $\Sigma_0=[E_0,\infty)$ for some $E_0\in\bbR$ and suppose that
$H_0$ has spectrum
\begin{equation}
\spec (H_0)=\Sigma_0. \lb{1.10}
\end{equation}
Then 
\begin{equation}
\cQ_0(x)=E_0\,\cI_m \text{ for a.e. $x\in\bbR$}. \lb{1.11}
\end{equation}
$(ii)$ Let $\Sigma_1=[E_0,E_1]\cup [E_2\infty)$ for some $E_0<E_1<E_2$
and suppose that $H_1$ has spectrum
\begin{equation}
\spec (H_1)=\Sigma_1. \lb{1.12}
\end{equation}
Then 
\begin{align}
\cQ_1(x)&=(1/3)(E_0+E_1+E_2) \cI_m \no \\
& \quad+2\cU
\diag(\wp(x+\omega_3+\alpha_1),\dots,\wp(x+\omega_3+\alpha_m))\cU^{-1} 
\lb{1.13} \\ 
& \quad \text{ for some $\alpha_j\in\bbR$, $1\leq j\leq m$ and a.e.
$x\in\bbR$,} \no
\end{align}
where $\cU$ is an $m\times m$ unitary matrix independent of $x\in\bbR$.
In particular, $\cQ_1$ satisfies the stationary KdV equation
\begin{equation}
\cQ'''_1-3(\cQ_1^2)'+2(E_0+E_1+E_2)\cQ'_1=0. \lb{1.14}
\end{equation}
\end{theorem}

As shown in \cite{CGHL00}, periodic Schr\"odinger operators in 
$L^2(\bbR)^{m\times m}$ with spectra of uniform (maximal) multiplicity
$2m$ are reflectionless in the sense that \eqref{1.4} holds for all
$\lambda$ in the open interior of the spectrum. Hence one obtains the
following result.

\begin{theorem} \lb{t1.5a} Let $m\in\bbN$, suppose 
$\cQ_\ell=\cQ_\ell^* \in L^1_{\loc}(\bbR)^{m\times m}$ is periodic and
define the self-adjoint Schr\"odinger operators $H_\ell$, $\ell=0,1$ in
$L^2(\bbR)^{m\times m}$ as in \eqref{1.9}. Assume
\begin{equation}
\spec (H_\ell)=\Sigma_\ell, \quad \ell=0,1, \lb{1.15}
\end{equation}
and suppose that $H_\ell$, $\ell=0,1$ has uniform $($maximal\,$)$
spectral multiplicity $2m$. Then $\cQ_\ell$, $\ell=0,1,$ are
reflectionless and hence the assertions \eqref{1.11}, \eqref{1.13}, and
\eqref{1.14} of Theorem \ref{t1.5} hold.
\end{theorem}

\begin{remark} \lb{rem1.6} ${}$ \\
$(i)$ The assumption of uniform $($maximal\,$)$ spectral multiplicity $2m$
in Theorem \ref{t1.5a}\,$(i)$ is an essential one. Otherwise, one can
easily  construct nonconstant potentials $\cQ$ such that the associated
Schr\"odinger operator $H_\cQ$ has overlapping band spectra and hence
spectrum equal to a half-line. For such a construction it suffices to
consider the case in which $\cQ$ is a diagonal matrix. In the special
scalar case
$m=1$,  reflectionless potentials automatically give rise to maximum
uniform spectral multiplicity $2$ for the corresponding scalar
Schr\"odinger operator in $L^2(\bbR)$. \\
$(ii)$ Theorem \ref{t1.5a}\,$(i)$ assuming $\cQ_0\in
L^{\infty}(\bbR)^{m \times m}$ to be periodic has  been proved by
Depr\'es \cite{De95} using an entirely different approach based on a 
detailed Floquet analysis. Depr\'es' result was reproved in
\cite{CGHL00} under the current general assumptions on $\cQ_0$
using methods based on matrix-valued Herglotz functions and trace
formulas. \\
$(iii)$ For different proofs of Borg's Theorem~\ref{t1.1} in the scalar
case $m=1$ we refer to \cite{Ho65}, \cite{Jo82}, \cite{Jo87}, 
\cite{Ko84}. \\
$(iv)$ Without loss of generality we focus on the limit point case of 
the differential expression $-\cI_m\f{d^2}{dx^2}+\cQ$ at $\pm \infty$ in
this paper. In fact, by a result originally due to Povzner \cite{Po53}, 
scalar  Schr\"odinger differential expressions leading to minimal
operators bounded from below are in the limit point case at $\pm\infty$.
Povzner's result was later also proved by Wienholtz \cite{Wi58} and is
reproduced as Theorem\ 35 in \cite[p.\ 58]{Gl65}. As shown in
\cite{CG02a}, Wienholtz's proof extends to the matrix case at hand. Since
the spectra $\Sigma_\ell$, $\ell=0,1$ are bounded from below, the limit
point assumption is justified $($and natural\,$)$. 
\end{remark}

This paper is another modest contribution to the inverse spectral theory
of matrix-valued Schr\"odinger (and Dirac-type) operators and part of a
recent program in this area (cf.\ \cite{CG01}, \cite{CG02}, \cite{CG02a},
\cite{CGHL00}, \cite{GH97}, \cite{GKM01}, \cite{GS02}, \cite{GS00}, and
\cite{GT00}). For other relevent recent literature in this context we
refer, for instance, to \cite{Ca00}, \cite{Ca00a}, \cite{Ca01},
\cite{Ca01a}, \cite{Ch98}, \cite{CLW01}, \cite{CS97}, \cite{JL98},
\cite{JL99}, \cite{JNO00}, \cite{Ma99}, \cite{Ma99a}, \cite{RK01},
\cite{Sh00}, \cite{Sh01}, \cite{SS98}. For the applicability of this
circle of ideas to the nonabelian Korteweg--de Vries hierarchy we refer
to \cite{GS02} and the references cited therein.

In Section \ref{s2} we recall basic facts on Weyl--Titchmarsh theory and
pencils of matrices as needed in the remainder of this paper. Section
\ref{s3} summarizes the principal results of paper \cite{GS02}.
Finally, in Section \ref{s4} we present the matrix extensions of Borg's
and Hochstadt's theorem and the corresponding connection with the
stationary KdV hierarchy. 

\section{Preliminaries on Weyl--Titchmarsh Theory \\ and Pencils of 
Matrices} \lb{s2}

The basic assumption for this section will be the following.

\begin{hypothesis}\lb{h2.1} Fix $m\in\bbN$, suppose $\cQ=\cQ^*\in
L_{\loc}^1(\bbR)^{m\times m}$, introduce the differential expression
\begin{equation}
\cL=-\cI_m\f{d^2}{dx^2}+\cQ, \quad x\in\bbR. \lb{2.1}
\end{equation} 
and suppose $\cL$ is in the limit point case at $\pm\infty$. 
\end{hypothesis}

Given Hypothesis~\ref{h2.1} we consider the matrix-valued
Schr\"odinger equation 
\begin{equation}
-\psi^{\prime\prime}(z,x)+\cQ(x)\psi(z,x) =z\psi(z,x)  
\text{ for a.e. $x\in\bbR$},  \lb{2.6}
\end{equation}
where $z\in\bbC$ plays the role of a spectral parameter and $\psi$ is
assumed to satisfy
\begin{equation}
\psi(z,\cdot), \psi^\prime (z,\cdot) \in \AC_{\loc}(\bbR)^{m\times m}.
\lb{2.7}
\end{equation}
Throughout this paper, $x$-derivatives are abbreviated by a prime
$\prime$.

Let $\Psi(z,x,x_0)$ be a $2m\times 2m$ normalized fundamental system of 
solutions of \eqref{2.6} at some $x_0\in\bbR$ which we partition as
\begin{equation}
\Psi(z,x,x_0)
=\begin{pmatrix}\theta(z,x,x_0) & \phi(z,x,x_0)\\
\theta'(z,x,x_0)& \phi'(z,x,x_0)\end{pmatrix}. \lb{2.16}
\end{equation}
Here $\prime$ denotes $d/dx$, $\theta(z,x,x_0)$ and $\phi(z,x,x_0)$ are
$m\times m$ matrices, entire with respect to $z\in\bbC$, and normalized
according to $\Psi(z,x_0,x_0)=\cI_{2m}$. 

By Hypothesis~\ref{h2.1}, the $m\times m$ Weyl--Titchmarsh matrices
associated with $\cL$, the half-lines $[x,\pm\infty)$, and a Dirichlet
boundary condition at $x$, are given by (c.f.\ \cite{HS81}, \cite{HS83},
\cite{HS84}, \cite{KR74}, \cite{Or76}, \cite{Sa92}, \cite{Sa94a},
\cite{Sa99a})
\begin{equation}
\cM_{\pm}(z,x)=\Psi_{\pm}'(z,x,x_0)\Psi_\pm(z,x,x_0)^{-1}, 
\quad z\in\bbC\backslash\bbR, \lb{2.47}
\end{equation}
where $\Psi_{\pm}(z,\cdot,x_0)$ satisfy
$(\cL-z\cI_m)\Psi_{\pm}(z,\cdot,x_0)=0$ and 
\begin{equation}
\Psi_\pm(z,\cdot,x_0)\in L^2([x_0,\pm\infty))^{m\times m}. \lb{2.20g}
\end{equation}

Next, we recall the definition of matrix-valued Herglotz function.

\begin{definition} \lb{d2.2}
A map $\cM\colon\bbC_+\to \bbC^{n\times n}$, $n\in\bbN$, extended to
$\bbC_-$ by $\cM(\bar z)= \cM(z)^*$ for all $z\in\bbC_+$, is called an
$n\times n$ Herglotz matrix if it is analytic on $\bbC_+$ and
$\Im(\cM(z))\ge 0$ for all $z\in\bbC_+$.
\end{definition}

\noindent Here we denote $\Im(\cM)=(\cM-\cM^*)/2i$ and
$\Re(\cM)=(\cM+\cM^*)/2$.

$\pm \cM_\pm(\cdot,x)$ are $m\times m$ Herglotz matrices of rank $m$ and
hence admit the representations
\begin{equation}
\pm \cM_{\pm}(z,x)=\Re(\pm \cM_{\pm}(\pm i,x)) +\int_\bbR 
d\Omega_\pm(\lambda,x) \,
\big((\lambda-z)^{-1}-\lambda(1+\lambda^2)^{-1}\big), \lb{2.25a} 
\end{equation}
where
\begin{equation}
\int_\bbR \|d\Omega_\pm(\lambda,x)\|_{\bbC^{m\times m}}\,
(1+\lambda^2)^{-1}<\infty \lb{2.27} 
\end{equation}
and
\begin{equation}
\Omega_\pm((\lambda,\mu],x)=\lim_{\delta\downarrow
0}\lim_{\varepsilon\downarrow 0}\f1\pi
\int_{\lambda+\delta}^{\mu+\delta} d\nu \, \Im(\pm
\cM_\pm(\nu+i\varepsilon,x)). \lb{2.29} 
\end{equation}

Necessary and sufficient conditions for $\cM_\pm(\cdot,x_0)$ to be the
half-line $m\times m$ Weyl--Titchmarsh matrix associated with a
Schr\"odinger operator on $[x_0,\pm\infty)$ in terms of the
corresponding measures $\Omega_\pm(\cdot,x_0)$ in the Herglotz
representation \eqref{2.25a} of $\cM_\pm(\cdot,x_0)$ can be derived using
the matrix-valued extension of the classical inverse spectral theory
approach due to Gelfand and Levitan \cite{GL51}, as worked out by
Rofe-Beketov \cite{Ro60}. The following result describes sufficient
conditions for a monotonically nondecreasing matrix function to be the
matrix spectral function of a half-line Schr\"odinger operator. It extends
well-known results in the scalar case $m=1$ (cf.\
\cite[Sects.\ 2.5, 2.9]{Le87}, \cite{LG64}, \cite[Sect.\ 26.5]{Na68},
\cite{Th79}).

\begin{theorem} [\cite{Ro60}] \lb{t2.3} 
Suppose $\Omega_+(\cdot,x_0)$ is a monotonically nondecreasing $m\times m$
matrix-valued function on $\bbR$. Then $\Omega_+(\cdot,x_0)$ is the
matrix  spectral function of a self-adjoint Schr\"odinger operator $H_+$
in $L^2([x_0,\infty))^{m}$ associated with the $m\times m$ matrix-valued
differential expression $\cL_+=-d^2/dx^2 \cI_m+\cQ$, $x>x_0$, with a
Dirichlet boundary condition at $x_0$, a self-adjoint boundary condition
at $\infty$ $($if necessary$)$, and a self-adjoint potential matrix $\cQ$
with $\cQ^{(r)}\in L^1([x_0,R])^{m\times m}$ for all $R>x_0$ if and
only if the following two conditions hold. \\
\noindent $(i)$ 
Whenever $f\in C([x_0,\infty))^{m\times 1}$ with compact support 
contained in $[x_0,\infty)$ and
\begin{equation}
\int_\bbR F(\lambda)^*d\Omega_+(\lambda,x_0)\,F(\lambda) =0, 
 \text{ then $f=0$~a.e.,} \lb{2.35}
\end{equation} 
where 
\begin{equation}
F(\lambda)=\lim_{R\uparrow\infty}\int_{x_0}^R
dx\,\f{\sin(\lambda^{1/2}(x-x_0))}{\lambda^{1/2}}f(x), \quad
\lambda\in\bbR.
\lb{2.36}
\end{equation}
\noindent $(ii)$ Define 
\begin{equation}
\wti\Omega_+(\lambda,x_0)=\begin{cases}
\Omega_+(\lambda,x_0)-\f{2}{3\pi}\lambda^{3/2}, & \lambda\geq 0 \\
\Omega_+(\lambda,x_0), & \lambda<0 \end{cases} \lb{2.37}
\end{equation}
and assume the limit 
\begin{equation}
\lim_{R\uparrow\infty}\int_{-\infty}^R d\wti\Omega_+(\lambda,x_0) \,
\f{\sin(\lambda^{1/2}(x-x_0))}{\lambda^{1/2}}= \Phi(x) \lb{2.38}
\end{equation}
exists and $\Phi\in L^\infty([x_0,R])^{m\times m}$ for all $R>x_0$.
Moreover, suppose that for some $r\in\bbN_0$, $\Phi^{(r+1)}\in
L^1([x_0,R])^{m\times m}$ for all $R>x_0$, and $\Phi(x_0)=0$.  
\end{theorem}

Assuming Hypothesis \ref{h2.1}, we next introduce the self-adjoint 
Schr\"odinger operator $H$ in $L^2(\bbR)^m$ by
\begin{align}
&H=-\cI_m \f{d^2}{dx^2}+\cQ, \lb{2.39} \\
&\dom(H)=\{g\in L^2(\bbR)^m \mid g,g^\prime\in
\AC_{\loc}(\bbR)^m;\, (-g^{\prime\prime}+\cQ g)\in L^2(\bbR)^m\}. \no
\end{align}
The resolvent of $H$ then reads 
\begin{equation}
((H-z)^{-1}f)(x)= \int_\bbR dx^\prime\, \cG(z,x,x^\prime)f(x^\prime), 
\quad z\in\bbC\backslash\bbR, \; f\in L^2(\bbR)^m, \lb{2.38a} 
\end{equation}
with the Green's matrix $\cG(z,x,x')$ of $H$ given by 
\begin{align}
\cG(z,x,x^\prime)=\Psi_\mp(z,x,x_0)[\cM_-(z,x_0) &
-\cM_+(z,x_0)]^{-1}\Psi_\pm(\overline z,x^\prime,x_0)^*, 
\no \\
& \hspace*{2.05cm} x\lesseqgtr x^\prime,\; 
z\in\bbC\backslash\bbR. \lb{2.33}
\end{align}
Introducing
\begin{equation}
\cN_\pm (z,x)=\cM_-(z,x)\pm \cM_+(z,x),
\quad z\in\bbC\backslash\bbR, \; x\in\bbR, \lb{2.33A} 
\end{equation}
the $2m\times 2m$ Weyl--Titchmarsh function $\cM(z,x)$ associated with
$H$ on $\bbR$ is then given by 
\begin{align}
\cM(z,x)&=\big(\cM_{p,q}(z,x)\big)_{p,q=1,2} \no \\
&=\begin{pmatrix}
\cM_\pm (z,x)\cN_-(z,x)^{-1}\cM_\mp (z,x)   
&\cN_-(z,x)^{-1}\cN_+(z,x)/2  \\ 
\cN_+(z,x)\cN_-(z,x)^{-1}/2  &
\cN_-(z,x)^{-1} \end{pmatrix}, \lb{2.34} \\
& \hspace*{7cm} z\in\bbC\backslash\bbR, \; x\in\bbR. \no
\end{align} 

Then $\cM(z,x)$ is a $2m\times 2m$ matrix-valued Herglotz function of
rank $2m$ with representations 
\begin{align}
\cM(z,x)&=\Re(\cM(i,x)) +\int_\bbR d\Omega(\lambda,x)\,
\big((\lambda-z)^{-1}-\lambda(1+\lambda^2)^{-1}\big) \lb{2.42} \\
&=\exp\bigg(\cC(x)+\int_\bbR d\lambda\,\Upsilon(\lambda,x)
\big((\lambda-z)^{-1}-\lambda(1+\lambda^2)^{-1}\big)\bigg), \lb{2.43}
\end{align}
where
\begin{align}
&\int_\bbR \Vert d\Omega(\lambda,x)
\Vert_{\bbC^{2m\times 2m}} \,(1+\lambda^2)^{-1}<\infty, \lb{2.44} \\
& \cC(x)=\cC(x)^*, \quad 0\leq \Upsilon(\cdot,x)\leq \cI_{2m} \text{ a.e.}
\end{align}
and
\begin{align}
\Omega((\lambda,\mu],x)&=\lim_{\delta\downarrow
0}\lim_{\varepsilon\downarrow 0}\f1\pi
\int_{\lambda+\delta}^{\mu+\delta} d\nu \, 
\Im(\cM(\nu+i\varepsilon,x)), \lb{2.46} \\
\Upsilon(\lambda,x)&=\lim_{\varepsilon\downarrow 0} \pi^{-1}
\Im(\ln(\cM(\lambda+i\varepsilon,x))) \text{ for a.e. $\lambda\in\bbR$.} 
\lb{2.46A}
\end{align}

The Herglotz, and particularly exponential Herglotz property (cf.\
\cite{AD56}, \cite{Ca76}, \cite{GT00}) of the diagonal Green's function 
of $H$, 
\begin{equation}
\gg(z,x)=\cG(z,x,x), \quad z\in\bbC\backslash\spec(H), \; x\in\bbR, 
\lb{2.46B}
\end{equation}
will be of particular importance in Section \ref{s4} and hence we note for
subsequent purpose,
\begin{equation}
\gg(z,x)=\exp\bigg(\gC(x)+\int_\bbR d\lambda\,\Xi(\lambda,x)
\big((\lambda-z)^{-1}-\lambda(1+\lambda^2)^{-1}\big)\bigg), \lb{2.47a} 
\end{equation}
where
\begin{align}
&\gC(x)=\gC(x)^*, \quad 0\leq \Xi(\cdot,x)\leq \cI_{m} \text{ a.e.,} 
\lb{2.47b} \\
& \Xi(\lambda,x)=\lim_{\varepsilon\downarrow 0} \pi^{-1}
\Im(\ln(\gg(\lambda+i\varepsilon,x))) \text{ for a.e. $\lambda\in\bbR$.} 
\lb{2.47c}
\end{align}

We also recall the following characterization of $\cM(z,x_0)$ to be 
used later. In the scalar context $m=1$, this characterization has been
used by Rofe-Beketov \cite{Ro67}, \cite{Ro91} (see also 
\cite[Sect.\ 7.3]{Le87}).

\begin{theorem} [\cite{Ro67}, \cite{Ro91}] \lb{t2.4} 
Assume Hypothesis~\ref{h2.1}, suppose that $z\in\bbC 
\backslash \bbR$, $x_0\in\bbR$, and let $\ell, r\in\bbN_0$. Then the
following assertions are equivalent. \\
$(i)$ $\cM(z,x_0)$ is the $2m\times 2m$ 
Weyl--Titchmarsh matrix associated with a Schr\"odinger operator 
$H$ in $L^2(\bbR)^m$ of the type \eqref{2.39} with an $m\times m$
matrix-valued potential $\cQ\in L^1_{\loc}(\bbR)$ and $\cQ\in
C^\ell((-\infty,x_0))$ and $\cQ\in C^r((x_0,\infty))$. \\
$(ii)$ $\cM(z,x_0)$ is of the type \eqref{2.34} with $\cM_\pm(z,x_0)$
being half-line $m\times m$ Weyl--Titchmarsh matrices on $[x_0,\pm\infty)$
corresponding to a Dirichlet boundary condition at $x_0$ and a
self-adjoint boundary condition at $-\infty$ and/or $\infty$ $($if
any$)$ which are associated with an $m\times m$ matrix-valued potential
$\cQ$ satisfying $\cQ\in C^\ell((-\infty,x_0))$ and $\cQ\in
C^r((x_0,\infty))$, respectively. \\ 
If $(i)$ or $(ii)$ holds, then the $2m\times 2m$ matrix-valued spectral
measure $\Omega(\cdot,x_0)$ associated with $\cM(z,x_0)$ is determined
by \eqref{2.34} and \eqref{2.46}.
\end{theorem}

Next, we consider variations of the reference point $x\in\bbR$. Since
$\Psi_\pm$ satisfies the second-order linear $m\times m$ matrix-valued
differential equation \eqref{2.6}, $\cM_\pm$ in \eqref{2.47} satisfies the
matrix-valued Riccati-type equation (independently of any limit point
assumptions at $\pm\infty$)
\begin{equation}
\cM_{\pm}'(z,x)+\cM_\pm(z,x)^2=\cQ(x)-z \cI_m, \quad x\in\bbR,
\;  z\in\bbC\backslash\bbR. \lb{2.49}
\end{equation}

The asymptotic high-energy behavior of $\cM_\pm(z,x)$ as $|z|\to\infty$
has recently been determined in \cite{CG01} under minimal smoothness
conditions on $\cQ$ and without assuming that $\cL$ is in the limit point
case at $\pm\infty$. Here we recall just a special case of the asymptotic
expansion proved in \cite{CG01} which is most suited for our discussion
at hand. We denote by $C_\varepsilon\subset\bbC_+$ the open sector with
vertex at zero, symmetry axis along the positive imaginary axis, and
opening angle $\varepsilon$, with $0<\varepsilon< \pi/2$.

\begin{theorem} [\cite{CG01}] \lb{t2.5}
Fix $x_0\in\bbR$ and let $x\geq x_0$. In addition to
Hypothesis~\ref{h2.1} suppose that 
$\cQ\in C^\infty([x_0,\pm\infty))^{m\times m}$ and that $\cL$ is in the
limit point case at $\pm\infty$. Let $\cM_\pm(z,x)$, $x\geq x_0$, be
defined as in \eqref{2.47}. Then, as $\abs{z}\to\infty$ in
$C_\varepsilon$, $\cM_\pm(z,x)$ has an asymptotic expansion of the form
$(\Im(z^{1/2})>0$, $z\in\bbC_+)$
\begin{equation}
\cM_\pm(z,x)\underset{\substack{\abs{z}\to\infty\\ z\in
C_\varepsilon}}{=} \pm i \cI_m z^{1/2}+\sum_{k=1}^N
\cM_{\pm,k}(x)z^{-k/2}+ o(|z|^{-N/2}), \quad N\in\bbN. \lb{2.50}
\end{equation}
The expansion \eqref{2.50} is uniform with respect to $\arg\,(z)$ for 
$|z|\to \infty$ in $C_\varepsilon$ and uniform in $x$ as long as $x$
varies in compact subsets of $[x_0,\infty)$. The expansion coefficients
$\cM_{\pm,k}(x)$ can be recursively computed from 
\begin{align}
\begin{split}
\cM_{\pm,1}(x)&=\mp\f{i}{2} \cQ(x),
\quad \cM_{\pm,2}(x)= \f1{4} \cQ^\prime(x),  \\
\cM_{\pm,k+1}(x)&=\pm\f{i}2\bigg(\cM_{\pm,k}^\prime(x)+
\sum_{\ell=1}^{k-1}\cM_{\pm,\ell}(x) \cM_{\pm,k-\ell}(x) \bigg),
\quad k\ge 2. \lb{2.51}
\end{split}
\end{align} 
The asymptotic expansion \eqref{2.50} can be differentiated to any
order with respect to $x$. \\
If one only assumes Hypothesis~\ref{h2.1} $($i.e., $\cQ\in
L^1 ([x_0,R])^{m\times m}$ for all $R>x_0$$)$, then 
\begin{equation}
\cM_\pm(z,x)\underset{\substack{\abs{z}\to\infty\\ z\in
C_\varepsilon}}{=} \pm i \cI_m z^{1/2}+ \oh(1). \lb{2.51a}
\end{equation}
\end{theorem}

\begin{remark} \lb{r2.6} 
Due to the recursion relation \eqref{2.51}, the coefficients
$\cM_{\pm,k}$ are universal polynomials in $\cQ$ and its $x$-derivatives
$($i.e., differential polynomials in $\cQ$$)$. 
\end{remark}

Finally, in addition to \eqref{2.33} (still assuming Hypothesis
\ref{h2.1}), one infers for the $2m\times 2m$  Weyl--Titchmarsh function
$\cM(z,x)$ associated with $H$ on $\bbR$ in connection with arbitrary
half-lines $[x,\pm\infty)$, $x\in\bbR$, 
\begin{align}
\cM(z,x)&=\big(\cM_{j,j^\prime}(z,x)\big)_{j,j^\prime=1,2},
\quad  z\in\bbC\backslash\bbR,  \lb{2.52} \\
\cM_{1,1}(z,x)&=\cM_\pm(z,x)[\cM_-(z,x)-\cM_+(z,x)]^{-1}
\cM_\mp(z,x) \no \\ 
&=\psi_+'(z,x,x_0)[\cM_-(z,x_0)-\cM_+(z,x_0)]^{-1}\psi_-'(\ol z,x,x_0)^*,
\lb{2.53} \\
\cM_{1,2}(z,x)&=2^{-1} [\cM_-(z,x)-\cM_+(z,x)]^{-1}
[\cM_-(z,x)+\cM_+(z,x)] \no \\
&=\psi_+(z,x,x_0)[\cM_-(z,x_0)-\cM_+(z,x_0)]^{-1}\psi_-'(\ol z,x,x_0)^*,
\lb{2.54} \\
\cM_{2,1}(z,x)&=2^{-1}
[\cM_-(z,x)+\cM_+(z,x)][\cM_-(z,x)-\cM_+(z,x)]^{-1} \no \\
&=\psi_+'(z,x,x_0)[\cM_-(z,x_0)-\cM_+(z,x_0)]^{-1}\psi_-(\ol z,x,x_0)^*,
\lb{2.55} \\
\cM_{2,2}(z,x)&=[\cM_-(z,x)-\cM_+(z,x)]^{-1} \no \\
&=\psi_+(z,x,x_0)[\cM_-(z,x_0)-\cM_+(z,x_0)]^{-1}\psi_-(\ol z,x,x_0)^*.
\lb{2.56} 
\end{align}
Introducing the convenient abbreviation,
\begin{equation}
\cM(z,x)=\begin{pmatrix} \gh(z,x) & -\gg_2(z,x) \\ -\gg_1(z,x) & \gg(z,x)
\end{pmatrix}, \quad z\in\bbC\backslash\bbR, \; x\in\bbR, \lb{2.57}
\end{equation}
one then verifies from \eqref{2.52}--\eqref{2.57} and from $\cM(\ol
z,x)^*=\cM(z,x)$, $\cM_\pm(\ol z,x)^*=\cM_\pm(z,x)$ that 
\begin{align}
&\gg(\ol z,x)^*=\gg(z,x), \quad \gg_2(\ol z,x)^*=\gg_1(z,x), 
\quad \gh(\ol z,x)^*=\gh(z,x), \lb{2.58} \\
&\gg(z,x)\gg_1(z,x)=\gg_2(z,x)\gg(z,x), \lb{2.59} \\
&\gh(z,x)\gg_2(z,x)=\gg_1(z,x)\gh(z,x), \lb{2.59a} \\
&\gg(z,x)=[\cM_-(z,x)-\cM_+(z,x)]^{-1}, \lb{2.60} \\
&\gg(z,x)\gh(z,x)-\gg_2(z,x)^2=-(1/4)\cI_m, \lb{2.61} \\
&\gh(z,x)\gg(z,x)-\gg_1(z,x)^2=-(1/4)\cI_m, \lb{2.62} \\ 
&\cM_\pm(z,x)=\mp(1/2)\gg(z,x)^{-1}-\gg(z,x)^{-1}\gg_2(z,x)  \lb{2.63} \\
&\qquad \qquad \;\;\;\,\, \,\mp(1/2)\gg(z,x)^{-1}-\gg_1(z,x)\gg(z,x)^{-1},
\lb{2.64}
\end{align}
assuming Hypothesis \ref{h2.1}. Moreover, \eqref{2.58}--\eqref{2.64} and
the Riccati-type equations \eqref{2.49} imply the following results for
$z\in\bbC\backslash\bbR$ and a.e.\ $x\in\bbR$,
\begin{align}
\gg'&=-(\gg_1+\gg_2), \lb{2.65} \\
\gg_1'&=-(\cQ-z\cI_m)\gg-\gh \lb{2.66} \\
& = (-\gg''+\gg\cQ-\cQ\gg)/2, \lb{2.66a} \\ 
\gg_2'&=-\gg(\cQ-z\cI_m)-\gh \lb{2.67} \\
& =(-\gg''+\cQ\gg-\gg\cQ)/2, \lb{2.67a} \\
\gh'&=-\gg_1(\cQ-z\cI_m)-(\cQ-z\cI_m)\gg_2, \lb{2.68} \\  
\gh&=[\gg''-\gg(\cQ-z\cI_m)-(\cQ-z\cI_m)\gg]/2 \lb{2.69} 
\end{align}
if $\cQ\in L^1_{\loc}(\bbR)^{m\times m}$, and 
\begin{align}
\gg_1''&=-2(\cQ-z\cI_m)\gg'-\cQ'\gg+\gg_1\cQ-\cQ\gg_1, \lb{2.66b} \\
\gg_2''&=-2\gg'(\cQ-z\cI_m)-\gg\cQ'+\cQ\gg_2-\gg_2\cQ \lb{2.67b}  
\end{align}
if in addition $\cQ'\in L^1_{\loc}(\bbR)^{m\times m}$. 

We conclude this section by recalling the definition of reflectionless
matrix-valued potentials as discussed in \cite{CGHL00}, \cite{GKM01},
\cite{GT00}, and \cite{KS88}. We follow the corresponding notion
introduced in connection with scalar Schr\"odinger operators and
refer to \cite{Cr89}, \cite{DS83}, \cite{KK88}, \cite{SY95} for
further details in this context.

\begin{definition} \lb{d2.7}
Assume Hypothesis \ref{h2.1} and define the self-adjoint Schr\"odinger
operator $H$ in $L^2(\bbR)^{m\times m}$ as in \eqref{2.39}. Suppose that 
$\spec_{\ess}(H)\neq\emptyset$ and let $\Xi$ be defined by \eqref{2.47c}.
Then $\cQ$ is called {\it reflectionless} if for all $x\in\bbR$, 
\begin{equation}
\Xi(\lambda,x) =(1/2)\cI_m \text{ for a.e.\
$\lambda\in\spec_{\ess}(H)$.} \lb{2.70}
\end{equation}
\end{definition}
\noindent Explicit examples of reflectionless potentials will be discussed
in Section \ref{s3}. If $\cQ$ is reflectionless we will sometimes slightly
abuse notation  and also call the corresponding Schr\"odinger operator
$H$ in $L^2(\bbR)^{m\times m}$ reflectionless.

\section{A Class of Matrix-Valued Schr\"odinger Operators \\
with Prescribed Finite-Band Spectra} \lb{s3}

Given the preliminaries of Section \ref{s2}, we now recall the
construction of a class of matrix-valued Schr\"odinger operators with a
prescribed finite-band spectrum of uniform maximum multiplicity, the
principal  result of \cite{GS02}. Let 
\begin{equation}
\{E_\ell\}_{0\leq \ell\leq 2n}\subseteq \bbR, \; n\in\bbN, \text{ with
$E_\ell<E_{\ell+1}$, $0\leq \ell\leq 2n-1$,} \lb{4.1}
\end{equation}
and introduce the polynomial 
\begin{equation}
R_{2n+1}(z)=\prod_{\ell=0}^{2n} (z-E_\ell), \quad z\in\bbC. \lb{4.2}
\end{equation}
Moreover, we define the square root of $R_{2n+1}$ by 
\begin{equation}
R_{2n+1}(\lambda)^{1/2}=\lim_{\varepsilon\downarrow 0}
R_{2n+1}(\lambda+i\varepsilon)^{1/2}, 
\quad \lambda\in\bbR, \lb{4.3}
\end{equation}
and
\begin{align}
R_{2n+1} (\lambda)^{1/2} &= |R_{2n+1} (\lambda)^{1/2} |\begin{cases}
(-1)^n i& \text{for $\lambda \in (-\infty, E_0)$},\\
(-1)^{n+j} i & \text{for $\lambda \in (E_{2j-1}, E_{2j})$},\; j=1,\dots,
n,\\ (-1)^{n+j}& \text{for $\lambda \in (E_{2j}, E_{2j+1})$}, \;
j=0,\dots, n-1,\\ 1 & \text{for $\lambda \in (E_{2n}, \infty)$},
\end{cases} \no \\
& \hspace*{8.5cm} \lambda\in\bbR \lb{4.4}
\end{align}
and analytically continue $R_{2n+1}^{1/2}$ {}from $\bbR$ to all of
$\bbC\backslash\Sigma_n$, where $\Sigma_n$ is defined by
\begin{equation}
\Sigma_n=\Bigg\{\bigcup_{j=0}^{n-1} [E_{2j},E_{2j+1}]\Bigg\}\cup
[E_{2n},\infty).  \lb{4.4a}
\end{equation}
In this context we also mention the useful formula
\begin{equation}
\ol{R_{2n+1}(\ol z)^{1/2}}=-R_{2n+1}(z)^{1/2}, \quad z\in\bbC_+.
\lb{4.4b} 
\end{equation}

\begin{theorem} [\cite{GS02}] \lb{t3.1}
Let $z\in\bbC\backslash\Sigma_n$ and $n\in\bbN$. Define $R_{2n+1}^{1/2}$
as in \eqref{4.1}--\eqref{4.4} followed by an analytic continuation to
$\bbC\backslash\Sigma_n$. Moreover, let $F_n$ and $H_{n+1}$ be two monic 
polynomials of degree $n$ and $n+1$, respectively. Then 
$iR_{2n+1}(z)^{-1/2}F_n(z)$ 
is a Herglotz function if and only if all zeros of $F_n$ are real and
there is precisely one zero in each of the intervals $[E_{2j-1},E_{2j}]$,
$1\leq j\leq n$. Moreover, if $iR_{2n+1}^{-1/2}F_n$ is a Herglotz
function, then it can be represented in the form
\begin{equation}
\f{iF_n(z)}{R_{2n+1}(z)^{1/2}}=\f{1}{\pi} \int_{\Sigma_n} 
\f{F_n(\lambda)d\lambda}{R_{2n+1}(\lambda)^{1/2}}\f{1}{\lambda-z}, 
\quad z\in\bbC\backslash\Sigma_n. \lb{4.5a}
\end{equation}
Similarly, $iR_{2n+1}(z)^{-1/2}H_{n+1}(z)$ 
is a Herglotz function if and only if all zeros of $H_{n+1}$ are real and
there is precisely one zero in each of the intervals
$(-\infty,E_0]$ and $[E_{2j-1},E_{2j}]$, $1\leq j\leq n$. Moreover, 
if $iR_{2n+1}^{-1/2}H_{n+1}$ is a Herglotz
function, then it can be represented in the form
\begin{align}
\f{iH_{n+1}(z)}{R_{2n+1}(z)^{1/2}}&=
\Re\bigg(\f{iH_{n+1}(i)}{R_{2n+1}(i)^{1/2}}\bigg) \no \\
& \quad + \f{1}{\pi} \int_{\Sigma_n} 
\f{H_{n+1}(\lambda)d\lambda}{R_{2n+1}(\lambda)^{1/2}}
\bigg(\f{1}{\lambda-z}-\f{\lambda}{1+\lambda^2}\bigg), 
\quad z\in\bbC\backslash\Sigma_n. \lb{4.6a}
\end{align}
\end{theorem}

Actually, Theorem \ref{t4.1} can be improved by invoking ideas developed
in the Appendix of \cite{KN77} (cf.\ also \cite{SY95}). Since this
appears to be of independent interest we provide a brief discussion. 

We start with the elementary observation that the Herglotz function
\begin{align}
m(z)&= \begin{cases}
\f{z-\beta}{z-\alpha}, &  -\infty <\alpha<\beta <\infty, \\
z-\beta, & \alpha=-\infty, \; \beta\in\bbR, \\
\f{-1}{z-\alpha}, & \alpha\in\bbR, \; \beta=+\infty, \end{cases} \\
& \hspace*{3.75cm} z\in\bbC_+, \no
\end{align}
admits the (exponential) representation
\begin{equation}
m(z)=C(\alpha,\beta)\exp{\int_\alpha^\beta
d\lambda\bigg(\f{1}{\lambda-z}-\f{\lambda}{1+\lambda^2}\bigg )},
\quad z\in \bbC_+, \lb{3.10}
\end{equation}
where
\begin{equation}
C(\alpha,\beta)= \begin{cases}
\Big(\f{1+\beta^2}{1+\alpha^2}\Big)^{1/2}, &  -\infty <\alpha<\beta
<\infty, \\
1, & \alpha=-\infty, \; \beta\in\bbR \text{ or }
\alpha\in\bbR, \; \beta=+\infty.
\end{cases}
\end{equation}

\begin{theorem} \lb{t3.1a}
Let $z\in\bbC\backslash\Sigma_n$, $n\in\bbN$, and define $R_{2n+1}^{1/2}$
as in \eqref{4.1}--\eqref{4.4} followed by an analytic continuation to
$\bbC\backslash\Sigma_n$. Suppose $M$ is a Herglotz function such that
\begin{equation}
\lim_{\varepsilon \downarrow 0} M(\lambda +i\varepsilon)\in i \bbR 
\; \text{ for a.e.\ $\lambda \in \Sigma_n$}
\end{equation}
and assume in addition that $M$ is real-valued on
$\bbC\backslash\Sigma$. Then $M$ is either of the form
\begin{equation}
M(z)=\f{i\hatt F_n(z)}{R_{2n+1}(z)^{1/2}}, \lb{4.5}
\end{equation}
where $\hatt F_n$ is  a polynomial of degree $n$ $($not necessarily
monic$)$, positive on the semi-axis
$(E_{2n}, \infty)$, with precisely one zero in each of the intervals
$[E_{2j-1},E_{2j}]$, $1\leq j\leq n$,  or else, $M$ is of the form
\begin{equation}
M(z)=\f{i\hatt H_{n+1}(z)}{R_{2n+1}(z)^{1/2}}, \lb{4.6}
\end{equation}
where $\hatt H_{n+1}$ is  a polynomial of degree $n+1$ $($not necessarily
monic$)$, positive on the semi-axis $(E_{2n}, \infty)$, with precisely
one zero in each of the intervals $(-\infty,E_0]$ and
$[E_{2j-1},E_{2j}]$, $1\leq j\leq n$. \\
\indent Moreover, if $i\hatt F_n/R_{2n+1}^{1/2}$ is a Herglotz
function, it can be represented in the form \eqref{4.5a}. 
Similarly, if $i\hatt H_{n+1}/R_{2n+1}^{1/2}$ is a Herglotz
function, it can be represented in the form \eqref{4.6a}.
\end{theorem}
\begin{proof} The Herglotz function $M$ admits the exponential
representation (cf.\ \cite{AD56})
\begin{equation}
M(z)=K\exp{\bigg
(\frac{1}{2}\int_{\Sigma_n}d\lambda\bigg(\f{1}{\lambda-z}
-\f{\lambda}{1+\lambda^2}\bigg)
+\int_{\Sigma_{n,-}}d\lambda\bigg(\f{1}{\lambda-z}
-\f{\lambda}{1+\lambda^2}\bigg)\bigg )},
\end{equation}
where $K>0$ and
\begin{equation}
\Sigma_{n,-}=\{ \lambda \in \bbR\backslash \Sigma_n \, | \,
M(\lambda)<0\}.
\end{equation}
Since the Herglotz function $M$ is strictly monotonically
increasing on $\bbR\backslash\Sigma_n$, $M$ can have at most one zero in
each interval $(-\infty,E_0)$, $(E_{2j-1},E_{2j})$, $j=1,\dots,n$.
Moreover, $M(E_0-0)$, $M(E_{2j}-0)\in\{+\infty,0\}$, 
$M(E_{2j-1}+0)\in\{-\infty,0\}$, $j=1,\dots,n$. Thus, depending on
whether or not $M\big|_{(-\infty,E_0)}\geq 0$, the set
$\Sigma_{n,-}$  admits one of the following two  representations
\begin{equation}\lb{sigma1}
\Sigma_{n,-}=\bigcup_{ j=1}^{n}(E_{2j-1}, \mu_{j})
\end{equation}
or
\begin{equation} \lb{sigma2}
\Sigma_{n,-}=(-\infty, \nu_0)\cup \bigcup_{ j=1}^{n}(E_{2j-1}, \nu_{j})
\end{equation}
for some $\nu_0\in (-\infty, E_0]$ and some 
$\mu_j,\nu_j\in [E_{2j-1},E_{2j}]$, $1\le j \le n$. Repeated use of
\eqref{3.10} then proves the representation
\begin{align}
M(z)&=C_1\bigg
(\prod_{j=0}^{n-1}\f{(E_{2j+1}-z)}{(E_{2j}-z)}\f{1}{(E_{2n}-z)}\bigg
)^{1/2}
\prod_{j=1}^n \f{(z-\mu_j)}{(z-E_{2j-1})} \no \\
&=C_1\f{\prod_{j=1}^n (z-\mu_j)}
{\Big (\prod_{\ell=0}^{2n}(z-E_\ell)\Big)^{1/2}} 
=\f{i\hatt F_n(z)}{R_{2n+1}(z)^{1/2}}, 
\end{align}
where
\begin{equation}
\hatt F_n(z)=C_1\prod_{j=1}^n (z-\mu_j)
\end{equation}
for some $C_1>0$, whenever \eqref{sigma1} holds, and
\begin{equation}
M(z)=C_2(z-\nu_0)\bigg
(\prod_{j=0}^{n-1}\f{E_{2j+1}-z}{E_{2j}-z}\f{1}{E_{2n}-z}\bigg)^{1/2}
\prod_{j=1}^n
\f{z-\nu_j}{z-E_{2j-1}}=\f{i\hatt H_{n+1}(z)}{R_{2n+1}(z)^{1/2}},
\end{equation}
where
\begin{equation}
\hatt H_{n+1}(z)=C_2\prod_{k=0}^n (z-\nu_k)
\end{equation}
for some $C_2>0$, whenever \eqref{sigma2} holds.
\end{proof}

Given $m\in\bbN$, we denote by
\begin{equation}
\cA(z)=\sum_{k=0}^n \cA_k z^k, \quad \cA_k\in\bbC^{m\times m}, 
\, 1\leq k\leq n, \; z\in\bbC, \lb{3.1}
\end{equation}
a polynomial pencil of $m\times m$ matrices (in short, a pencil) in the
following. $\cA$ is called  of degree $n\in\bbN_0$ if $\cA_n\neq 0$ and
monic if $\cA_n=\cI_m$. 

\begin{definition} \lb{d3.2} Let $A$ be a pencil of the type 
\eqref{3.1}. \\
$(i)$ The pencil $\cA$ is called {\it self-adjoint} if
$\cA_k=\cA_k^*$ for all $1\leq k\leq n$ $($i.e., $\cA(\ol z)^*=\cA(z)$ for
all $z\in\bbC$$)$. \\
$(ii)$ A self-adjoint pencil $\cA$ is called {\it weakly hyperbolic} if
$\cA_n>0$ and for all $f\in\bbC^m\backslash\{0\}$, the roots of the
polynomial $(f,\cA(\cdot)f)_{\bbC^m}$ are real. If in addition all these
zeros are distinct, the pencil $\cA$ is called {\it hyperbolic}. \\
$(iii)$ Let $\cA$ be a weakly hyperbolic pencil and denote by 
$\{p_j(\cA,f)\}_{1\leq j\leq n}$,
\begin{equation}
 p_j(\cA,f)\leq
p_{j+1}(\cA,f), \; 1\leq j\leq n-1, \; f\in\bbC^m\backslash\{0\},
\end{equation}
the roots of the polynomial $(f,\cA(\cdot)f)_{\bbC^m}$ ordered in
magnitude. The range of the roots 
$p_j(\cA,f)$, $f\in\bbC^m\backslash\{0\}$ is denoted by
$\Delta_j(\cA)$ and called the {\it $j$th root zone} of $\cA$. \\
$(iv)$ A hyperbolic pencil $\cA$ is called {\it strongly hyperbolic}
if $\ol {\Delta_j(\cA)}$ and $\ol {\Delta_k(\cA)}$ are mutually
disjoint for $j\neq k$, $1\leq j,k \leq n$. 
\end{definition}

For details on spectral theory of polynomial matrix (in fact, operator)
pencils we refer, for instance, to \cite{Ma88a}, \cite{MM76},
\cite{MM85}. 

\begin{corollary} [\cite{GS02}] \lb{c3.3} 
Let $z\in\bbC\backslash\Sigma_n$ and $m,n\in\bbN$. Define 
$R_{2n+1}^{1/2}$ as in \eqref{4.1}--\eqref{4.4} followed by an analytic
continuation to $\bbC\backslash\Sigma_n$. Moreover let $\cF_n$ and
$\cH_{n+1}$ be two monic $m\times m$ matrix pencils of degree $n$ and
$n+1$, respectively. Then $(i/2)R_{2n+1}^{-1/2}\cF_n$ is a Herglotz
matrix if and only if the root zones $\Delta_j(\cF_n)$ of $\cF_n$
satisfy
\begin{equation}
\Delta_j(\cF_n)\subseteq [E_{2j-1},E_{2j}], \quad 1\leq j \leq n. 
\lb{4.12a}
\end{equation}
Analogously, $(i/2)R_{2n+1}^{-1/2}\cH_{n+1}$ is a Herglotz matrix if
and only if the root zones $\Delta_j(\cH_{n+1})$ of $\cH_{n+1}$
satisfy
\begin{equation}
\Delta_0(\cH_{n+1})\subset (-\infty,E_{0}], \quad
\Delta_j(\cH_{n+1})\subseteq [E_{2j-1},E_{2j}], \quad 1\leq j \leq n. 
\lb{4.12b}
\end{equation}
If \eqref{4.12a} $($resp., \eqref{4.12b}$)$ holds, then $\cF_n$
$($resp., $\cH_{n+1}$$)$ is a strongly hyperbolic pencil.
\end{corollary}

Next, we define the following $2m\times 2m$ matrix $\cM_{\Sigma_n}(z,x_0)$
which will turn out to be the underlying Weyl--Titchmarsh matrix
associated with the class of $m \times m$ matrix-valued Schr\"odinger
operators with prescribed finite-band spectrum $\Sigma_n$ of maximal
multiplicity. We introduce, for fixed $x_0\in\bbR$,
\begin{align}
\cM_{\Sigma_n}(z,x_0)&=\big(\cM_{\Sigma_n,p,q}(z,x_0)
\big)_{1\leq p,q\leq 2} \lb{4.13} \\
&=\f{i}{2R_{2n+1}(z)^{1/2}}\begin{pmatrix} \cH_{n+1,\Sigma_n}(z,x_0)
& -\cG_{2,n-1,\Sigma_n}(z,x_0)  \\ 
-\cG_{1,n-1,\Sigma_n}(z,x_0) & \cF_{n,\Sigma_n}(z,x_0) \end{pmatrix}, 
\; z\in\bbC\backslash\Sigma_n. \no
\end{align}
Here $R_{2n+1}(z)^{1/2}$ is defined as in \eqref{4.1}--\eqref{4.4}
followed by analytic continuation into $\bbC\backslash\Sigma$  and the
polynomial matrix pencils $\cF_{n,\Sigma_n}$, $\cG_{1,n-1,\Sigma_n}$,
$\cG_{2,n-1,\Sigma_n}$, and $\cH_{n+1,\Sigma_n}$ are introduced as
follows: \\

\noindent $(i)$ $\cF_{n,\Sigma_n}(\cdot,x_0)$ is an $m\times m$ monic
matrix pencil of degree $n$, that is, $\cF_{n,\Sigma}(\cdot,x_0)$
is of the type 
\begin{equation}
\cF_{n,\Sigma_n}(z,x_0)=\sum_{\ell=0}^n \cF_{n-\ell,\Sigma_n}(x_0) z^\ell,
\quad \cF_{0,\Sigma_n}(x_0)=\cI_m, \; z\in\bbC \lb{4.19}
\end{equation}
and 
\begin{equation}
\f{i}{2R_{2n+1}^{1/2}}\cF_{n,\Sigma_n}(\cdot,x_0) \text{ is assumed to be
an
$m\times m$ Herglotz matrix.}  \lb{4.20}  
\end{equation}
Hence $\cF_{n,\Sigma_n}(\cdot,x_0)$ is a self-adjoint $($in fact,
strongly hyperbolic$)$ pencil,
\begin{equation}
\cF_{n,\Sigma_n}(\ol z,x_0)^*=\cF_{n,\Sigma_n}(z,x_0), \quad z\in\bbC
\lb{4.19a}
\end{equation}
and $(i/2)R_{2n+1}^{-1/2}\cF_{n,\Sigma_n}$ and
$2iR_{2n+1}^{1/2}\cF_{n,\Sigma_n}^{-1}$ admit the Herglotz representations
\begin{align}
&\f{i}{2R_{2n+1}(z)^{1/2}}\cF_{n,\Sigma_n}(z,x_0)
=\f{1}{2\pi} \int_{\Sigma_n} 
\f{d\lambda}{R_{2n+1}(\lambda)^{1/2}}\cF_{n,\Sigma_n}(\lambda,x_0) 
\f{1}{\lambda-z}, 
\quad z\in\bbC\backslash\Sigma_n, \lb{4.21} \\
&iR_{2n+1}(z)^{1/2}\cF_{n,\Sigma_n}(z,x_0)^{-1} \no \\
&=\f{1}{\pi} \int_{\Sigma_n} d\lambda\,
R_{2n+1}(\lambda)^{1/2}\cF_{n,\Sigma_n}(\lambda,x_0)^{-1} 
\bigg(\f{1}{\lambda -z}-\f{\lambda}{1+\lambda^2}\bigg) \no \\
&\quad +\Gamma_{\Sigma_n,0}(x_0)-\sum_{k=1}^{N} 
(z-\mu_k (x_0))^{-1} \Gamma_{\Sigma_n,k}(x_0) , \lb{4.22} \\
& \hspace*{2.35cm} z\in\bbC\backslash
\{\Sigma_n\cup\{\mu_k(x_0)\}_{1\leq k\leq N}\}, \no
\end{align}
where
\begin{align}
&\Gamma_{\Sigma_n,0}(x_0)=\Gamma_{\Sigma_n,0}(x_0)^*\in\bbC^{m\times m},
\quad 0\leq \Gamma_{\Sigma_n,k}(x_0)\in\bbC^{m\times m}, \; 1\leq k \leq
N, \no \\ 
& \sum_{k=1}^N \rank(\Gamma_{\Sigma_n,k}(x_0)) \leq mn, \quad 
\mu_k(x_0)\in \bigcup_{j=1}^n [E_{2j-1},E_{2j}], \quad 1\leq k \leq N.
\lb{4.26}
\end{align}
In fact, there are precisely $m$ numbers $\mu_k(x_0)$ in 
$[E_{2j-1},E_{2j}]$ for each $1 \leq j\leq n$, counting multiplicity (they
are the points $z$  where $\cF_{n,\Sigma_n}(z,x_0)$ is not invertible). \\

\noindent $(ii)$ Given these facts we now define  
\begin{align}
\cG_{1,n-1,\Sigma_n}(z,x_0)&=\bigg(\sum_{k=1}^N
\f{\varepsilon_k (x_0)}{z-\mu_k(x_0)}
\Gamma_{\Sigma_n,k}(x_0)\bigg)\cF_{n,\Sigma_n}(z,x_0),
\lb{4.27} \\
\cG_{2,n-1,\Sigma_n}(z,x_0)&=\cF_{n,\Sigma_n}(z,x_0) \bigg(\sum_{k=1}^N
\f{\varepsilon_k (x_0)}{z-\mu_k(x_0)} \Gamma_{\Sigma_n,k}(x_0)\bigg),
\lb{4.28} \\ 
& \hspace*{-1.6cm} \varepsilon_k (x_0)\in\{1,-1\}, \; 1\leq k \leq N, 
\quad z\in\bbC\backslash\{\mu_k(x_0)\}_{1\leq k\leq N}, \lb{4.29} 
\end{align}
and
\begin{align}
&\cH_{n+1,\Sigma_n}(z,x_0)=R_{2n+1}(z) \cF_{n,\Sigma_n}(z,x_0)^{-1}
\lb{4.30}  \\ 
& +\bigg(\sum_{k=1}^N
\f{\varepsilon_k (x_0)}{z-\mu_k(x_0)}
\Gamma_{\Sigma_n,k}(x_0)\bigg)\cF_{n,\Sigma_n}(z,x_0) 
\bigg(\sum_{\ell=1}^N
\f{\varepsilon_\ell (x_0)}{z-\mu_\ell(x_0)}
\Gamma_{\Sigma_n,\ell}(x_0)\bigg),  \no \\
& \hspace*{7.5cm} z\in\bbC\backslash\{\mu_k(x_0)\}_{1\leq k\leq N}. \no
\end{align}

\begin{lemma} [\cite{GS02}] \lb{l3.4} Let
$z\in\bbC\backslash\{\mu_k(x_0)\}_{1\leq k\leq N}$.
$\cG_{p,n-1,\Sigma_n}(\cdot,x_0)$, $p=1,2$, are $m\times m$ polynomial
matrix pencils of equal degree at most $n-1$ and
$\cH_{n+1,\Sigma_n}(\cdot,x_0)$ is a strongly hyperbolic $($and hence
self-adjoint$)$ $m\times m$ monic matrix pencil of degree $n+1$.
Moreover, the following identities hold.
\begin{align}
&\cG_{2,n-1,\Sigma_n}(\ol z,x_0)^*=\cG_{1,n-1,\Sigma_n}(z,x_0),
\lb{4.30a} \\
&\cF_{n,\Sigma_n}(z,x_0)\cG_{1,n-1,\Sigma_n}(z,x_0)
=\cG_{2,n-1,\Sigma_n}(z,x_0)\cF_{n,\Sigma_n}(z,x_0), \lb{4.31} \\
&\cH_{n+1,\Sigma_n}(z,x_0)\cG_{2,n-1,\Sigma_n}(z,x_0)
=\cG_{1,n-1,\Sigma_n}(z,x_0)\cH_{n+1,\Sigma_n}(z,x_0), \lb{4.31a} \\
&\cF_{n,\Sigma_n}(z,x_0)\cH_{n+1,\Sigma_n}(z,x_0)
-\cG_{2,n-1,\Sigma_n}(z,x_0)^2
=R_{2n+1}(z)\cI_m, \lb{4.32} \\ 
&\cH_{n+1,\Sigma_n}(z,x_0)\cF_{n,\Sigma_n}(z,x_0)
-\cG_{1,n-1,\Sigma_n}(z,x_0)^2 =R_{2n+1}(z)\cI_m. \lb{4.33} 
\end{align}
\end{lemma}

Next, introducing
\begin{subequations} \lb{4.40}
\begin{align}
&\cM_{\pm,\Sigma_n}(z,x_0) \no \\
&= \pm
iR_{2n+1}(z)^{1/2}\cF_{n,\Sigma_n}(z,x_0)^{-1}
-\cG_{1,n-1,\Sigma_n}(z,x_0)\cF_{n,\Sigma_n}(z,x_0)^{-1} \lb{4.40a} \\
&=\pm
iR_{2n+1}(z)^{1/2}\cF_{n,\Sigma_n}(z,x_0)^{-1}
-\cF_{n,\Sigma_n}(z,x_0)^{-1}\cG_{2,n-1,\Sigma_n}(z,x_0), 
\lb{4.40b} \\
& \hspace*{5.5cm} 
z\in\bbC\backslash\{\Sigma_n\cup\{\mu_k(x_0)\}_{1\leq k\leq N}\}, \no 
\end{align}
\end{subequations}
$\pm\cM_{\pm,\Sigma_n}(\cdot,x_0)$ are $m\times m$ Herglotz
matrices with representations
\begin{align}
\pm \cM_{\pm,\Sigma_n}(z,x_0)&= \f{1}{\pi} \int_{\Sigma_n} d\lambda\,
R_{2n+1}(\lambda)^{1/2}\cF_{n,\Sigma_n}(\lambda,x_0)^{-1} 
\bigg(\f{1}{\lambda -z}-\f{\lambda}{1+\lambda^2}\bigg) \no \\
& \quad +\Gamma_{\Sigma_n,0}(x_0) - \sum_{k=1}^N
\f{1 \pm\varepsilon_k (x_0)}{z-\mu_k(x_0)} \Gamma_{\Sigma_n,k}(x_0),
\lb{4.43} \\ 
& \hspace*{1.9cm} z\in\bbC\backslash\{\Sigma_n\cup
\{\mu_k(x_0)\}_{1\leq k\leq N}\}. \no 
\end{align} 
Applying Theorem \ref{t2.3}, $\cM_{\pm,\Sigma_n}(z,x_0)$ are seen to be 
the half-line Weyl--Titchmarsh matrices uniquely associated
with a potential $\cQ_{\Sigma_n}\in L^1_{\loc}(\bbR)^{m\times m}$. In
addition, we denote by $\psi_{\pm,\Sigma_n}(z,x,x_0)$ the Weyl solutions
\eqref{2.47}, \eqref{2.20g} associated  with $\cQ_{\Sigma_n}$, 
\begin{equation}
\psi_{\pm,\Sigma_n}(z,x,x_0)=\theta_{\Sigma_n}(z,x,x_0)
+\phi_{\Sigma_n}(z,x,x_0)
\cM_{\pm,\Sigma_n}(z,x_0), \quad  z\in\bbC\backslash\Sigma_n, \lb{4.52a}
\end{equation}
where, in obvious notation, $\theta_{\Sigma_n}(z,x,x_0)$,
$\phi_{\Sigma_n}(z,x,x_0)$ denote the fundamental system \eqref{2.16}
corresponding to $\cQ_{\Sigma_n}$. The $2m\times 2m$ Weyl--Titchmarsh
matrix associated with $\cQ_{\Sigma_n}$ on $\bbR$ is then given by
\begin{align}
\cM_{\Sigma_n}(z,x)&=\big(\cM_{\Sigma,p,q}(z,x)\big)_{1\leq p,q\leq 2} 
 \lb{4.53} \\
&=\f{i}{2R_{2n+1}(z)^{1/2}}\begin{pmatrix} \cH_{n+1,\Sigma_n}(z,x)
& -\cG_{2,n-1,\Sigma_n}(z,x)  \\ 
-\cG_{1,n-1,\Sigma_n}(z,x) & \cF_{n,\Sigma_n}(z,x) \end{pmatrix}, 
\quad z\in\bbC\backslash\Sigma,  \no
\end{align}
where we abbreviated
\begin{align}
\cF_{n,\Sigma_n}(z,x)&=\theta_{\Sigma_n}(z,x,x_0)
\cF_{n,\Sigma_n}(z,x_0)\theta_{\Sigma_n}(\ol z,x,x_0)^* \no \\
& \quad + \phi_{\Sigma_n}(z,x,x_0)\cH_{n+1,\Sigma_n}(z,x_0)
\phi_{\Sigma_n}(\ol z,x,x_0)^* \no \\
& \quad - \phi_{\Sigma_n}(z,x,x_0)\cG_{1,n-1,\Sigma_n}(z,x_0)
\theta_{\Sigma_n}(\ol z,x,x_0)^* \no \\
& \quad - \theta_{\Sigma_n}(z,x,x_0)\cG_{2,n-1,\Sigma_n}(z,x_0)
\phi_{\Sigma_n}(\ol z,x,x_0)^*, \lb{4.54} \\
\cG_{1,n-1,\Sigma_n}(z,x)&=
-\theta_{\Sigma_n}'(z,x,x_0)\cF_{n,\Sigma_n}(z,x_0)
\theta_{\Sigma_n}(\ol z,x,x_0)^*
\no \\ & \quad - \phi_{\Sigma_n}'(z,x,x_0)\cH_{n+1,\Sigma_n}(z,x_0)
\phi_{\Sigma_n}(\ol z,x,x_0)^* \no \\
& \quad + \phi_{\Sigma_n}'(z,x,x_0)\cG_{1,n-1,\Sigma_n}(z,x_0)
\theta_{\Sigma_n}(\ol z,x,x_0)^* \no \\
& \quad + \theta_{\Sigma_n}'(z,x,x_0)\cG_{2,n-1,\Sigma_n}(z,x_0)
\phi_{\Sigma_n}(\ol z,x,x_0)^*, \lb{4.55} \\
\cG_{2,n-1,\Sigma_n}(z,x)&=-\theta_{\Sigma_n}(z,x,x_0)
\cF_{n,\Sigma_n}(z,x_0)\theta_{\Sigma_n}'(\ol z,x,x_0)^* \no \\
& \quad - \phi_{\Sigma_n}(z,x,x_0)\cH_{n+1,\Sigma_n}(z,x_0)
\phi_{\Sigma_n}'(\ol z,x,x_0)^* \no \\
& \quad + \phi_{\Sigma_n}(z,x,x_0)\cG_{1,n-1,\Sigma_n}(z,x_0)
\theta_{\Sigma_n}'(\ol z,x,x_0)^* \no \\
& \quad + \theta_{\Sigma_n}(z,x,x_0)\cG_{2,n-1,\Sigma_n}(z,x_0)
\phi_{\Sigma_n}'(\ol z,x,x_0)^*, \lb{4.56} \\
\cH_{n+1,\Sigma_n}(z,x)&=\theta_{\Sigma_n}'(z,x,x_0)
\cF_{n,\Sigma_n}(z,x_0)\theta_{\Sigma_n}'(\ol z,x,x_0)^* \no \\
& \quad + \phi_{\Sigma_n}'(z,x,x_0)\cH_{n+1,\Sigma_n}(z,x_0)
\phi_{\Sigma_n}'(\ol z,x,x_0)^* \no \\
& \quad - \phi_{\Sigma_n}'(z,x,x_0)\cG_{1,n-1,\Sigma_n}(z,x_0)
\theta_{\Sigma_n}'(\ol z,x,x_0)^* \no \\
& \quad - \theta_{\Sigma_n}'(z,x,x_0)\cG_{2,n-1,\Sigma_n}(z,x_0)
\phi_{\Sigma_n}'(\ol z,x,x_0)^*. \lb{4.57}
\end{align}
Considerations of this type can be found in \cite[Sect.\ 8.2]{Le87} 
in the special scalar case $m=1$ and in the matrix context $m\in\bbN$ in 
\cite[Sect.\ 9.4]{Sa99a}.

One then infers
\begin{align}
\cF_{n,\Sigma_n}'&=-(\cG_{1,n-1,\Sigma_n}+\cG_{2,n-1,\Sigma_n}), \lb{4.58} \\
\cG_{1,n-1,\Sigma_n}'&=-(Q_{\Sigma_n}-z\cI_m)\cF_{n,\Sigma_n}
-\cH_{n+1,\Sigma_n} \lb{4.59} \\
& =(-\cF_n''+\cF_n\cQ_{\Sigma_n}-\cQ_{\Sigma_n}\cF_{n,\Sigma_n})/2, 
\lb{4.59a} \\
\cG_{1,n-1,\Sigma_n}''&=-2(\cQ_{\Sigma_n}-z\cI_m)\cF_{n,\Sigma_n}'
-\cQ'_{\Sigma_n}\cF_{n,\Sigma_n}+\cG_{1,n-1,\Sigma_n}\cQ_{\Sigma_n}
-\cQ_{\Sigma_n}\cG_{1,n-1,\Sigma_n}, \lb{4.59b} \\
\cG_{2,n-1,\Sigma_n}'&=-\cF_{n,\Sigma_n}(Q_{\Sigma_n}-z\cI_m)
-\cH_{n+1,\Sigma_n} \lb{4.60} \\
&=(-\cF_{n,\Sigma_n}''+\cQ_{\Sigma_n}\cF_{n,\Sigma_n}
-\cF_{n,\Sigma_n}\cQ_{\Sigma_n})/2, \lb{4.60a} \\
\cG_{2,n-1,\Sigma_n}''&=-2\cF_{n,\Sigma_n}'(\cQ_{\Sigma_n}
-z\cI_m)-\cF_{n,\Sigma_n}\cQ'_{\Sigma_n}
+\cQ_{\Sigma_n}\cG_{2,n-1,\Sigma_n}-\cG_{2,n-1,\Sigma_n}\cQ_{\Sigma_n},
\lb{4.60b} \\ 
\cH_{n+1,\Sigma_n}'&=-\cG_{1,n-1,\Sigma_n}(Q_{\Sigma_n}-z\cI_m)
-(Q_{\Sigma_n}-z\cI_m)\cG_{2,n-1,\Sigma_n}, \lb{4.61} \\
\cH_{n+1,\Sigma_n}&=[\cF_{n,\Sigma_n}''
-\cF_{n,\Sigma_n}(\cQ_{\Sigma_n}-z\cI_m)-(\cQ_{\Sigma_n}
-z\cI_m)\cF_{n\Sigma_n}]/2 \lb{4.61a}
\end{align} 
and
\begin{align}
&\cF_{n,\Sigma_n}(\ol z,x)^*=\cF_{n,\Sigma_n}(z,x), \quad  
\cH_{n+1,\Sigma_n}(\ol z,x)^*=\cH_{n+1,\Sigma_n}(z,x), \no \\
&\cG_{2,n-1,\Sigma}(\ol z,x)^*=\cG_{1,n-1,\Sigma_n}(z,x), \lb{4.62} \\
&\cF_{n,\Sigma_n}(z,x)\cG_{1,n-1,\Sigma_n}(z,x)
=\cG_{2,n-1,\Sigma_n}(z,x)\cF_{n,\Sigma_n}(z,x), \lb{4.63} \\
&\cH_{n+1,\Sigma_n}(z,x)\cG_{2,n-1,\Sigma_n}(z,x)
=\cG_{1,n-1,\Sigma_n}(z,x)\cH_{n+1,\Sigma_n}(z,x), \lb{4.64} \\
&\cH_{n+1,\Sigma_n}(z,x)\cF_{n,\Sigma_n}(z,x)
-\cG_{1,n-1,\Sigma_n}(z,x)^2=R_{2n+1}(z)\cI_m, \lb{4.65} \\
&\cF_{n,\Sigma_n}(z,x)\cH_{n+1,\Sigma_n}(z,x)
-\cG_{2,n-1,\Sigma_n}(z,x)^2=R_{2n+1}(z)\cI_m. \lb{4.66}
\end{align}

Combining \eqref{2.52}--\eqref{2.56} and \eqref{4.53} then yields
\begin{subequations} \lb{4.67}
\begin{align}
&\cM_{\pm,\Sigma_n}(z,x) \no \\
&= \pm iR_{2n+1}(z)^{1/2}\cF_{n,\Sigma_n}(z,x)^{-1}
-\cG_{1,n-1,\Sigma_n}(z,x)\cF_{n,\Sigma_n}(z,x)^{-1} \lb{4.67a} \\
&=\pm iR_{2n+1}(z)^{1/2}\cF_{n,\Sigma_n}(z,x)^{-1}
-\cF_{n,\Sigma_n}(z,x)^{-1}\cG_{2,n-1,\Sigma_n}(z,x), \lb{4.67b} \\
& \hspace*{7.95cm} \quad z\in\bbC\backslash\bbR. \no 
\end{align}
\end{subequations}
One observes that for each $x\in\bbR$, $\cM_{+,\Sigma_n}(\cdot,x)$ is the
analytic continuation of $\cM_{-,\Sigma_n}(\cdot,x)$ through the set
$\Sigma_n$, and vice versa, 
\begin{align}
&\lim_{\varepsilon\downarrow 0}\cM_{+,\Sigma_n}(\lambda+i\varepsilon,x)
=\lim_{\varepsilon \downarrow 0}\cM_{-,\Sigma_n}(\lambda-i\varepsilon,x)
=\lim_{\varepsilon \downarrow
0}\cM_{-,\Sigma_n}(\lambda+i\varepsilon,x)^*,
\lb{4.67c} \\
& \hspace*{4.85cm} \lambda\in \bigcup_{j=0}^{n-1} (E_{2j},E_{2j+1}) \cup
(E_{2n},\infty), \; x\in\bbR. \no
\end{align}
In other words, for each $x\in\bbR$, $\cM_{+,\Sigma_n}(\cdot,x)$ and 
$\cM_{-,\Sigma_n}(\cdot,x)$ are the two branches of an analytic
matrix-valued function $\cM_{\Sigma_n}(\cdot,x)$ on the two-sheeted
Riemann surface of $R_{2n+1}^{1/2}$. This implies that $\cQ_{\Sigma_n}$
is reflectionless as will be discussed in Lemma \ref{l3.6}. In
addition, it is worthwhile to emphasize that in the present case of
reflectionless potentials, $\cF_{n,\Sigma_n}(z,x_0)$ and
$\cG_{1,n-1,\Sigma_n}(z,x_0)$ for some fixed $x_0\in\bbR$, uniquely
determine $\cM_{\pm,\Sigma_n}(z,x_0)$ and hence $\cQ_{\Sigma_n}(x)$ for
all $x\in\bbR$. 

Introducing the open interior $\Sigma^o_n$ of $\Sigma_n$ defined by
$\Sigma^o_n=\bigcup_{j=0}^{n-1} (E_{2j},E_{2j+1})\cup (E_{2n},\infty)$, 
one obtains the following results.

\begin{theorem} [\cite{GS02}] \lb{t3.5}
Let $z\in\bbC\backslash\bbR$ and $x\in\bbR$. Then \\
$(i)$ $\cF_{n,\Sigma_n}(\cdot,x)$ and $\cH_{n+1,\Sigma_n}(\cdot,x)$ are
strongly hyperbolic $($and hence self-adjoint$)$ $m\times m$ monic matrix
pencils of degree $n$ and
$n+1$, respectively, and $\cG_{p,n-1,\Sigma_n}(\cdot,x)$, $p=1,2$, are
$m\times m$ matrix pencils of degree $n-1$. \\
$(ii)$ The differential expression $\cL_{\Sigma_n}=-\cI_m\f{d^2}{dx^2} +
\cQ_{\Sigma_n}$ is in the limit point case at $\pm\infty$. \\
$(iii)$ $\cM_{\pm,\Sigma_n}(z,\cdot)$ in \eqref{4.67} satisfy the
matrix-valued Riccati-type equation
\begin{equation}
\cM_{\pm,\Sigma_n}'(z,x)+\cM_{\pm,\Sigma_n}(z,x)^2=\cQ_{\Sigma_n}(x)-z
\cI_m, \quad x\in\bbR, \;  z\in\bbC\backslash\bbR. \lb{4.68}
\end{equation}
Moreover, $\cM_{\pm,\Sigma_n}(z,x)$ in \eqref{4.67} are the
$m\times m$ Weyl--Titchmarsh matrices associated with self-adjoint
operators $H^D_{\pm,x,\Sigma_n}$ in $L^2([x,\pm\infty))^m$, with a
Dirichlet boundary condition at the point $x$ and an $m\times m$
matrix-valued potential
$\cQ_{\Sigma_n}$ satisfying 
\begin{equation}
\cQ_{\Sigma_n}=\cQ_{\Sigma_n}^*\in C^\infty(\bbR)^{m\times m}, \quad 
\cQ^{(r)}_{\Sigma_n}\in L^\infty(\bbR) \text{ for all $r\in\bbN_0$.}
\lb{4.43a}
\end{equation}
In addition, $\cQ_{\Sigma_n}$ is analytic in a neighborhood of the real
axis. $H^D_{\pm,x,\Sigma_n}$ is given by 
\begin{align}
&H_{\pm,x,\Sigma_n}^D=-\cI_m\f{d^2}{dx^2} + \cQ_{\Sigma_n}, \lb{4.44} \\
&\dom(H_{\pm,x,\Sigma_n}^D)=\{g\in L^2((x,\pm\infty))^m \,|\,g,g'\in
AC([x,c])^m \text{ for all } c\gtrless x;  \no \\
& \hspace*{3.2cm} \lim_{\varepsilon\downarrow 0} g(x\pm\varepsilon)=0;
\;(-g''+\cQ_{\Sigma_n} g)\in L^2((x,\pm\infty))^m \}. \no
\end{align}
$(iv)$ For each $x\in\bbR$, $\cM_{\Sigma_n}(z,x)$ in \eqref{4.53}
is a $2m\times 2m$ Weyl--Titchmarsh matrix associated with
the self-adjoint operator $H_{\Sigma_n}$ in $L^2(\bbR)^m$  defined by
\begin{align}
&H_{\Sigma_n}=-\cI_m \f{d^2}{dx^2}+\cQ_{\Sigma_n}, \lb{4.45} \\
&\dom(H_{\Sigma_n})=\{g\in L^2(\bbR)^m \mid g,g^\prime\in
\AC_{\loc}(\bbR)^m;\, (-g^{\prime\prime}+\cQ_{\Sigma_n} g)\in
L^2(\bbR)^m\}. \no
\end{align} 
In particular, $\cM_{\Sigma_n}(\cdot,x)$ is a $2m\times 2m$
Herglotz matrix of $H_{\Sigma_n}$ admitting a representation of the type
\eqref{2.42}, with measure $\Omega_{\Sigma_n}(\cdot,x)$ given by 
\begin{equation} 
d\Omega_{\Sigma_n}(\lambda,x)=\begin{cases} \f{1}{2\pi
R_{2n+1}(\lambda)^{1/2}}\left(\begin{smallmatrix}
\cH_{n+1,\Sigma_n}(\lambda,x) & -\cG_{2,n-1,\Sigma_n}(\lambda,x) \\
-\cG_{1,n-1,\Sigma_n}(\lambda,x) & \cF_{n,\Sigma_n}(\lambda,x) 
\end{smallmatrix}\right)d\lambda, & \lambda\in\Sigma^o_n, \\
0, & \lambda\in\bbR\backslash\Sigma_n. \end{cases} \lb{4.69}
\end{equation}
$(v)$ $H_{\Sigma_n}$ has purely absolutely continuous spectrum $\Sigma_n$,
\begin{equation}
\spec(H_{\Sigma_n})=\spec_{\text{ac}}(H_{\Sigma_n})=\Sigma_n, \quad 
\spec_{\p}(H_{\Sigma_n})=\spec_{\singc}(H_{\Sigma_n})=\emptyset, \lb{4.46}
\end{equation}
with $\spec(H_{\Sigma_n})$ of uniform spectral multiplicity $2m$.
\end{theorem}

It should be emphasized that the construction of $\cQ_{\Sigma_n}$ in the
scalar case $m=1$ is due to Levitan \cite{Le77} (see also \cite{Le77a},
\cite{Le84}, \cite[Ch.\ 8]{Le87}, \cite{LS88}).

That $\cQ_{\Sigma_n}$ is reflectionless is an elementary consequence of
\eqref{4.67c} as discussed next. 

\begin{lemma} [\cite{GS02}] \lb{l3.6}
Denote by $\gg_{\Sigma_n}(z,x)=\cG_{\Sigma_n}(z,x,x)$, $z\in\bbC_+$,
$x\in\bbR$, the diagonal Green's function of $H_{\Sigma_n}$. Then
\begin{equation}
\lim_{\varepsilon\downarrow 0}\gg_{\Sigma_n}(\lambda+i\varepsilon,x)=-
\lim_{\varepsilon\downarrow 0}\gg_{\Sigma_n}(\lambda+i\varepsilon,x)^*
\text{ for all $\lambda\in\Sigma_n^o$} \lb{3.59}
\end{equation}
and hence $\cQ_{\Sigma_n}$ is reflectionless. 
\end{lemma} 
\begin{proof}
Since $\gg(z,x)=(\cM_-(z,x)-\cM_+(z,x))^{-1}$, \eqref{4.67c} implies
\eqref{3.59}. The latter implies that 
$\lim_{\varepsilon\downarrow
0}\gg(\lambda+i\varepsilon,x)=iG_{\Sigma_n}(\lambda,x)$ for all
$\lambda\in\Sigma_n^o$ for some 
$m\times m$ matrix $G_{\Sigma_n}(\lambda,x)>0$. This in turn implies
\begin{equation}
\Xi_{\Sigma_n}(\lambda,x)=\lim_{\varepsilon\downarrow 0} \pi^{-1}
\Im(\ln(iG_{\Sigma_n}(\lambda+i\varepsilon,x)))=(1/2)\cI_m \text{ for all
$\lambda\in\Sigma^o$}, \lb{3.60}
\end{equation}
and hence $\cQ_{\Sigma_n}$ is reflectionless by Definition \ref{d2.7}.
\end{proof}

Next, we briefly turn to the stationary matrix Korteweg--de
Vries (KdV) hierarchy (cf.\ \cite[Ch.~15]{Di91}, \cite{GD77}) and show
that the finite-band potential
$\cQ_\Sigma$ satisfies some (and hence infinitely many) equations of the
stationary KdV equations. 

Assuming $\cQ=\cQ^*\in C^\infty(\bbR)^{m\times m}$, we recall the
expansions (cf.\ Theorem \ref{t2.5})
\begin{equation}
\gg(z,x)=[\cM_-(z,x)-\cM_+(z,x)]^{-1}\underset{\substack{\abs{z}\to\infty\\
z\in C_\varepsilon}}{=}\f{i}{2z^{1/2}}\sum_{k=0}^\infty
\hatt\gR_{k}(x)z^{-k} \lb{3.61} 
\end{equation}
for some coeficients $\hatt\gR_k$. Explicitly, one obtains
\begin{equation}
\hatt\gR_0=\cI_m, \quad \hatt\gR_1=\tfrac{1}{2}\cQ, \quad 
\hatt\gR_2=-\tfrac{1}{8}\cQ''+\tfrac{3}{8}\cQ^2, \text{ etc.} \lb{3.62}
\end{equation}
 
The stationary KdV hierarchy is then given by
\begin{equation}
\skdv_{k}(\cQ)=-2\sum_{\ell=0}^k c_{k-\ell}
\hatt\gR_{\ell+1}'(\cQ,\dots)=0, \quad k\in\bbN_0, \lb{3.63}
\end{equation}
where $\{c_\ell\}_{\ell=1,\dots,k}\subset\bbC$, $c_0=1$ denotes a set of
constants.

\noindent By Remark \ref{r2.6}, each $\hatt\gR_\ell$ is a differential
polynomial  in $\cQ$ and next we slightly abuse notation and indicate
this by writing $\hatt\gR_\ell(\cQ,\dots)$ for
$\hatt\gR_\ell(x)$, $\hatt\gR_{\ell+1}'(\cQ,\dots)$ for
$\hatt\gR_{\ell+1}'(x)$, etc. 

\begin{theorem} [\cite{GS02}] \lb{t3.8}
The self-adjoint finite-band potential $\cQ_{\Sigma_n}\in
C^\infty(\bbR)^{m\times m}$, discussed in Theorem \ref{t3.5}, is a
stationary KdV solution satisfying
\begin{equation}
\skdv_{n}(\cQ_{\Sigma_n})=-2\sum_{\ell=0}^n c_{n-\ell}(\ul E)
\hatt\gR_{\ell+1}'(\cQ_{\Sigma_n},\dots)=0. \lb{3.64}
\end{equation}
 
\end{theorem}  
Here $c_\ell(\ul E)$ are given by 
\begin{align}
c_0(\ul E)&=1,\no \\
c_k(\ul E)&=-\!\!\!\!\!\sum_{\substack{j_0,\dots,j_{2n}=0\\
 j_0+\cdots+j_{2n}=k}}^{k}\!\!\!\!\!
\f{(2j_0)!\cdots(2j_{2n})!}
{2^{2k} (j_0!)^2\cdots (j_{2n}!)^2 (2j_0-1)\cdots(2j_{2n}-1)}
E_0^{j_0}\cdots E_{2n}^{j_{2n}}, \no \\
& \hspace*{7cm} k=1,\dots,n. \label{3.65} 
\end{align}

\section{Matrix Extensions of Borg's and Hochstadt's Theorems} \lb{s4}

In this our principal section, we now prove Theorem \ref{t1.5}, the matrix
extension of Borg's and Hochstadt's theorem, Theorems \ref{t1.1} and
\ref{t1.3}. Our strategy of proof will be the following: First we show
that the (reflectionless) Schr\"odinger operators $H_{\Sigma_\ell}$
constructed in our previous Section \ref{s3} with spectrum $\Sigma_\ell$, 
satisfy the conclusions \eqref{1.11} and \eqref{1.13} for $\ell=0,1$,
respectively. Then, in a second step, we will prove that any
reflectionless Schr\"odinger operator with spectrum given by
$\Sigma_\ell$, $\ell=0,1$, is precisely of the form $H_{\Sigma_\ell}$ as
constructed in Section \ref{s3}.

\begin{theorem} \lb{t4.1}
Let $\ell=0,1$ and $\cQ_{\Sigma_\ell}$ be the finite-band potentials
constructed in Section \ref{s3}, with
$\spec(H_{\Sigma_\ell})=\Sigma_{\ell}$ $($cf.\ Theorem \ref{t3.5}$)$. 
Then
\begin{equation}
\cQ_{\Sigma_0}(x)=E_0\cI_m \text{ for a.e. $x\in\bbR$} \lb{5.1}
\end{equation}
and  
\begin{align}
\cQ_{\Sigma_1}(x)&=(1/3)(E_0+E_1+E_2) \cI_m \no \\
& \quad +2\cU
\diag(\wp(x+\omega_3+\alpha_1),\dots,\wp(x+\omega_3+\alpha_m))\cU^{-1} 
\lb{5.2} \\ 
& \quad \text{ for some $\alpha_j\in\bbR$, $1\leq j\leq m$ and a.e.
$x\in\bbR$,} \no
\end{align}
where $\cU$ is an $m\times m$ unitary matrix independent of $x\in\bbR$.
Moreover, $\cQ_{\Sigma_0}$ satisfies the first element of the KdV
hierarchy,
\begin{equation}
\cQ_{\Sigma_0}'=0, \lb{5.3}
\end{equation}
and $\cQ_{\Sigma_1}$ satisfies the stationary KdV equation
\begin{equation}
\cQ_{\Sigma_1}'''-3(\cQ_{\Sigma_1}^2)'
+2(E_0+E_1+E_2)\cQ_{\Sigma_1}'=0. \lb{5.4}
\end{equation}
\end{theorem}
\begin{proof}
We consider the elementary case $\ell=0$ first. Then the explicit
expressions,
\begin{align}
&\cF_{0,\Sigma_0}(z,x)=\cI_m, \quad \cG_{p,-1,\Sigma_0}(z,x)=0, \; p=1,2,
\quad \cH_{1,\Sigma_0}(z,x)=(z-E_0)\cI_m, \lb{5.5} \\
& \cM_{\pm,\Sigma_0}(z,x)=\pm i(z-E_0)^{1/2}\cI_m, \quad 
\gg_{\Sigma_0}(z,x)=(i/2)(z-E_0)^{-1/2}\cI_m, \text{ etc.,} \lb{5.6}
\end{align}
immediately imply \eqref{5.1} and \eqref{5.3}. Hence we turn to the 
case $\ell=1$. In this case one obtains,
\begin{align}
&\cF_{1,\Sigma_1}(z,x)=z\cI_m+(1/2)\cQ_{\Sigma_1}(x)+c_1\cI_m, \lb{5.7} \\
&\cG_{p,0,\Sigma_1}(z,x)=-(1/4)\cQ_{\Sigma_1}'(x), \; p=1,2, \lb{5.8} \\
&\cH_{2,\Sigma_1}(z,x)=z^2\cI_m+z(-(1/2)\cQ_{\Sigma_1}+c_1\cI_m)
+(1/4)\cQ_{\Sigma_1}''-(1/2)\cQ_{\Sigma_1}^2-c_1\cQ_{\Sigma_1}, \lb{5.9}
\\ 
& \cM_{\pm,\Sigma_1}(z,x)=\pm iR_3(z)^{1/2}\cF_{1,\Sigma_1}(z,x)^{-1}
-\cG_{1,0,\Sigma_1}(z,x)\cF_{1,\Sigma_1}(z,x)^{-1}, \lb{5.10} \\  
&\gg_{\Sigma_1}(z,x)=(i/2)R_3(z)^{-1/2}
(z\cI_m+(1/2)\cQ_{\Sigma_1}(x)+c_1\cI_m)^{-1}, \text{ etc.,} \lb{5.11}
\end{align} 
abbreviating
\begin{equation}
c_1=-(1/2)(E_0+E_1+E_2). \lb{5.12}
\end{equation}
Combining \eqref{4.63}, \eqref{5.7}, and \eqref{5.8}, $\cQ_{\Sigma_1}(x)$ 
and $\cQ_{\Sigma_1}'(x)$ commute and hence one obtains for each
$x\in\bbR$,
\begin{equation}
[\cQ_{\Sigma_1}^{(r)}(x),\cQ_{\Sigma_1}^{(s)}(x)]=0 \text{ for all
$r,s\in\bbN_0$}. \lb{5.13}
\end{equation}
Since $\cQ_{\Sigma_1}$ and all its derivatives are self-adjoint, one can
simultaneously diagonalize the family of matrices 
$\{\cQ_{\Sigma_1}^{(r)}(x_0)\}_{r\in\bbN_0}$ by a fixed unitary $m\times
m$ matrix $\cU$. By \eqref{5.7}--\eqref{5.10}, this also shows that 
$\cF_{1,\Sigma_1}(z,x_0)$, $\cG_{p,0,\Sigma_1}(z,x_0)$, $p=1,2$,
$\cH_{2,\Sigma_1}(z,x_0)$, and $\cM_{\pm,\Sigma_1}(z,x_0)$ can all be
simultaneously diagonalized by $\cU$. In particular, the spectral
measure in the Herglotz representation \eqref{4.43} of
$\cM_{\pm,\Sigma_1}(z,x_0)$ can be diagonalized by $\cU$. After
diagonalization with $\cU$, the inverse spectral approach in \cite{Ro60}
(i.e., the matrix-valued extension of the scalar Gelfand--Levitan method
\cite{GL51}, \cite[]{Le87}, \cite{Th79}) then yields a diagonal matrix
potential of the type 
\begin{equation}
\diag(q_{1,\Sigma_1}(x),\dots,q_{m,\Sigma_1}(x)), \lb{5.14}
\end{equation}
and hence $\cQ_{\Sigma_1}$ itself is of the form 
\begin{equation}
\cQ_{\Sigma_1}(x)=\cU\diag(q_{1,\Sigma_1}(x),
\dots,q_{m,\Sigma_1}(x))\cU^{-1}. \lb{5.15}
\end{equation}
In order to determine the scalar potentials $q_{k,\Sigma_1}$, 
$1\leq k\leq m$, it now suffices to solve the corresponding scalar
problem ($m=1$). But then Hochstadt's result \cite{Ho65} immediately
yields
\begin{equation}
q_{k,\Sigma_1}(x)=(1/3)(E_0+E_1+E_2)+2\wp(x+\omega_3+\alpha_k), \quad 
1\leq k\leq m, \lb{5.16}
\end{equation} 
for some $\{\alpha_k\}_{1\leq k\leq m}\subset\bbR$, and hence
\eqref{5.2}. Finally, \eqref{5.4} is a consequence of \eqref{3.64},
\eqref{3.65}, taking $n=1$. 
\end{proof}

To complete the proof of Theorem \ref{t1.5}, we next will prove that
reflectionless Schr\"odinger operators with spectrum equal to
$\Sigma_\ell$, $\ell=0,1$, in fact, coincide with some element of the
family $H_{\Sigma_\ell}$, described in Section \ref{s3}. 

\begin{theorem} \lb{t4.2}
Suppose $\cQ_\ell$, $\ell=0,1$ satisfies 
Hypothesis \ref{h2.1}, define $H_\ell$ as in \eqref{2.39}, and
suppose $\spec(H_\ell)=\Sigma_\ell$, $\ell=0,1$. In addition, assume that
$Q_\ell$, $\ell=0,1$, is reflectionless. Then
$H_\ell$ coincides with an element of the family of operators
$H_{\Sigma_\ell}$ parametrized by a choice of
$\cF_{\ell,\Sigma_\ell}(z,x_0)$, $\cG_{1,\ell-1,\Sigma_\ell}(z,x_0)$, 
$\ell=0,1$\footnote{More precisely, a choice of
$\cF_{0,\Sigma_\ell}(z,x_0)=\cI_m$ for $\ell=0$ and a choice of
$\cF_{1,\Sigma_\ell}(z,x_0)$ and a set of signs
$\varepsilon_k(x_0)\in\{1,-1\}$, $k=1,\dots,N$ (cf.\ \eqref{5.51b})
for $\ell=1$.}. 
\end{theorem}
\begin{proof}
We start with the case $\ell=0$. Recalling $\Sigma_0=[E_0,\infty)$, 
and $R_1(z)^{1/2}=(z-E_0)^{1/2}$ defined as in \eqref{4.3},\eqref{4.4},
followed by an analytic continuation from $\bbR$ to
$\bbC\backslash\Sigma_0$, we denote by $\gg(\Sigma_0,z,x)$ the diagonal
Green's function of $H_0$ (cf.\ \eqref{2.46B}),
\begin{equation}
\gg(\Sigma_0,z,x)=(\cM_{-}(\Sigma_0,z,x)-\cM_{+}(\Sigma_0,z,x))^{-1}, 
\lb{5.17}
\end{equation} 
where in obvious notation $\cM_\pm(\Sigma_0,z,x)$ denote the
Weyl-Titchmarsh matrices associated with $\cQ_0$. Since
$\cQ_0$ is reflectionless, the matrix
$\Xi(\Sigma_0,\cdot,x)$ in its associated exponential Herglotz
representation \eqref{2.47a} satisfies,
\begin{equation}
\Xi(\Sigma_0,\lambda,x)=\begin{cases}(1/2)\cI_m &\text{for a.e.\
$\lambda \in (E_0,\infty)$,} \\
0 &\text{for a.e.\ $\lambda
\in (-\infty,E_0)$.}\end{cases} \lb{5.18}
\end{equation}
Insertion of \eqref{5.18} into \eqref{2.47a} then yields 
\begin{equation}
\gg(\Sigma_0,z,x)=i(z-E_0)^{-1/2}\exp(\gC_0(x)), \quad
z\in\bbC\backslash\Sigma_0. \lb{5.19}
\end{equation}
A comparison with the high-energy asymptotics of $\gg$ implied by
\eqref{2.51a} and \eqref{5.17} yields 
\begin{equation}
\gg(\Sigma_0,z,x)\underset{\substack{\abs{z}\to\infty\\ z\in
C_\varepsilon}}{=} (i/2) \cI_m z^{-1/2}+\oh(1) \lb{5.20}
\end{equation}
and hence $\gC_0(x)=-\ln(2)\cI_m$ implying
\begin{equation}
\gg(\Sigma_0,z,x)=(i/2)(z-E_0)^{-1/2}\cI_m, \quad
z\in\bbC\backslash\Sigma_0. \lb{5.21}
\end{equation}
Thus,
\begin{equation}
\cQ_0(x)=\cQ_{\Sigma_0}(x)=E_0\,\cI_m, \quad x\in\bbR, \lb{5.22}
\end{equation}
with $\cQ_{\Sigma_0}$ constructed in Section \ref{s3} (cf.\ also 
\eqref{5.5}, \eqref{5.6}). \\ 
Next we turn to the case $\ell=1$. Recalling
$\Sigma_1=[E_0,E_1]\cup[E_2,\infty)$,  and
$R_3(z)^{1/2}=[(z-E_0)(z-E_1)(z-E_2)]^{1/2}$ defined as in
\eqref{4.3}, \eqref{4.4}, followed by an analytic continuation from
$\bbR$ to $\bbC\backslash\Sigma_1$, we introduce 
\begin{equation}
\gg(\Sigma_1,z,x)=(\cM_{-}(\Sigma_1,z,x)-\cM_{+}(\Sigma_1,z,x))^{-1}, 
\lb{5.23}
\end{equation} 
where in obvious notation $\cM_\pm(\Sigma_1,z,x)$ denote the
Weyl-Titchmarsh matrices associated with $\cQ_1$, and note
that $\gg(\Sigma_1,\lambda,x)$ is a self-adjoint matrix for 
$\lambda\in (-\infty,E_0)\cup (E_1,E_2)$, that is, for all $x\in\bbR$,
\begin{equation}
\gg(\Sigma_1,\lambda,x)=\gg(\Sigma_1,\lambda,x)^*, \; \lambda \in
(-\infty,E_0)\cup (E_1,E_2). \lb{5.24}
\end{equation}
Since by hypothesis $\cQ_1$ is reflectionless, one also has
for all $x\in\bbR$,
\begin{equation}
\lim_{\varepsilon\downarrow 0}\gg(\Sigma_1,\lambda+i\varepsilon,x)^*=
-\lim_{\varepsilon\downarrow 0}\gg(\Sigma_1,\lambda+i\varepsilon,x) 
\text{ for a.e.\ $\lambda\in (E_0,E_1)\cup (E_2,\infty)$.} \lb{5.25} 
\end{equation}
Combining \eqref{5.24} and \eqref{5.25} with the properties of
$R_3(z)^{1/2}$ as discussed in \eqref{4.4}, one infers that for all
$x\in\bbR$,  
\begin{equation}
-i\lim_{\varepsilon\downarrow 0}
R_3(\lambda+i\varepsilon)^{1/2}\gg(\Sigma_1,\lambda+i\varepsilon,x)
\text{ is self-adjoint for a.e.\ $\lambda\in\bbR$.} \lb{5.26}
\end{equation}
Next, we will take a closer look at $\gg(\Sigma_1,\cdot)$ and show that
\eqref{5.26} in fact holds for all $\lambda\in\bbR$. Since
$\cQ_1$ is reflectionless, the matrix
$\Xi(\Sigma_1,\cdot,x)$ in its associated exponential Herglotz
representation \eqref{2.47a} satisfies,
\begin{equation}
\Xi(\Sigma_1,\lambda,x)=\begin{cases}(1/2)\cI_m &\text{for a.e.\
$\lambda \in (E_0,E_1)\cup (E_2,\infty)$,} \\
0 &\text{for a.e.\ $\lambda \in (-\infty,E_0)$.}\end{cases}
\lb{5.27}
\end{equation}
Insertion of \eqref{5.27} into \eqref{2.47a} then yields 
\begin{align}
&\gg(\Sigma_1,z,x)=i\bigg(\f{1+E_0^2}{1+E_1^2}\bigg)^{1/4}\bigg(\f{(z-E_1)}{(z-E_0)(z-E_2)}\bigg)^{1/2}
\lb{5.28} \\ 
& \times 
\exp\bigg(\gC_1(x)+\int_{E_1}^{E_2}d\lambda\,\Xi(\Sigma_1,\lambda,x)
\big((\lambda-z)^{-1}-\lambda(1+\lambda^2)^{-1}\big)\bigg),
\quad z\in\bbC\backslash\Sigma_1. \no
\end{align}
As a consequence, one obtains
\begin{align}
&-iR_3(z)^{1/2}\gg(\Sigma_1,z,x)=\bigg(\f{1+E_0^2}{1+E_1^2}\bigg)^{1/4}
(z-E_1) \lb{5.29} \\
& \times 
\exp\bigg(\gC_1(x)+\int_{E_1}^{E_2}d\lambda\,\Xi(\Sigma_1,\lambda,x)
\big((\lambda-z)^{-1}-\lambda(1+\lambda^2)^{-1}\big)\bigg),
\quad z\in\bbC\backslash\Sigma_1. \no
\end{align}
By \eqref{5.29}, $-i\lim_{\varepsilon\downarrow 0}
R_3(\lambda+i\varepsilon)^{1/2}\gg(\Sigma_1,\lambda+i\varepsilon,x)$ is
self-adjoint for $\lambda\in\bbR\backslash (E_1,E_2)$. However, since
$\gg(\Sigma_1,z,x)$ is analytic in $z\in\bbC\backslash\Sigma_1$, one
infers that $-i\lim_{\varepsilon\downarrow 0}
R_3(\lambda+i\varepsilon)^{1/2}\gg(\Sigma_1,\lambda+i\varepsilon,x)$ is
self-adjoint for $\lambda\in\bbR\backslash\{E_1,E_2\}$. Next, a
comparison with the high-energy asymptotics of $\gg$ implied by
\eqref{2.51a} and \eqref{5.23} yields 
\begin{equation}
\gg(\Sigma_1,z,x)\underset{\substack{\abs{z}\to\infty\\ z\in
C_\varepsilon}}{=} (i/2) \cI_m z^{-1/2}+\oh(1) \lb{5.30}
\end{equation}
and hence 
\begin{equation}
\gC_1(x)=\bigg(\f{1}{4}\ln\bigg(\f{1+E_1^2}{1+E_0^2}\bigg)-\ln(2)\bigg)
\cI_m+\int_{E_1}^{E_2}
d\lambda\,\Xi(\Sigma_1,\lambda,x)\f{\lambda}{1+\lambda^2}. \lb{5.31}
\end{equation}
Moreover, since $0\leq\Xi(\Sigma_1,z,x)\leq\cI_m$, we obtain for
$z=E_1-\varepsilon$,
\begin{align} 
&(E_1-z)\exp\bigg(\gC_1(x)+\int_{E_1}^{E_2}d\lambda\,
\Xi(\Sigma_1,\lambda,x)\big((\lambda-z)^{-1}-\lambda(1+\lambda^2)^{-1}
\big)\bigg) \no\\
&=\f{E_2-E_1+\varepsilon}{2}\bigg(\f{(1+E_2^2)^2}{(1+E_0^2)(1+E_1^2)}
\bigg)^{1/4} \no \\
&\quad \times \exp\bigg(-\int_{E_1}^{E_2} d\lambda\,
(\cI_m-\Xi(\Sigma_1,\lambda,x))(\lambda-E_1+\varepsilon)^{-1}\bigg) 
\lb{5.32}
\end{align}
and hence \eqref{5.29} remains bounded at $z=E_1$. Thus,
$-i\lim_{\varepsilon\downarrow 0}
R_3(\lambda+i\varepsilon)^{1/2}\gg(\Sigma_1,\lambda+i\varepsilon,x)$ is
self-adjoint for $\lambda\in\bbR\backslash\{E_2\}$ and by the
Schwartz reflection principle, $-i R_3(z)^{1/2}\gg(\Sigma_1,z,x)$ is
analytic for $z\in\bbC\backslash\{E_2\}$. Finally, if $E_2$ would be a
pole of $-i R_3(z)^{1/2}\gg(\Sigma_1,z,x)$, then $\gg(\Sigma_1,z,x)$
would have a $(z-E_2)^{-(3/2)}$ singularity at $E_2$, contradicting the
Herglotz property of $\gg(\Sigma_1,\cdot,x)$. (Of course, the same
argument applies to $z=E_1$.) Thus, 
\begin{equation}
-i R_3(z)^{1/2}\gg(\Sigma_1,z,x) \text{ is entire with respect to $z$.}
\lb{5.33}
\end{equation}
By \eqref{5.29}, one infers the bound 
\begin{equation}
\|-i R_3(z)^{1/2}\gg(\Sigma_1,z,x)\|\leq C(x)|z| \text{ for
$|z|>\max(|E_1|,|E_2|)$} \lb{5.34}
\end{equation}
for some constant $C(x)>0$. Thus, $\gg(\Sigma_1,z,x)$ is of the form
\begin{equation}
\gg(\Sigma_1,\lambda,x)=(i/2)R_3(z)^{-1/2}(\cA(x)z+\cB(x)), \lb{5.35}
\end{equation}
for some $m\times m$ matrices $\cA(x),\cB(x)\in\bbC^{m\times m}$. Hence
$\gg(\Sigma_1,z,x)$ has an asymptotic expansion to all orders as
$|z|\to\infty$ and an insertion of the asymptotic expansion into the
Riccati-type equation \eqref{2.49} yields (cf.\ also \eqref{2.51}) 
\begin{align}
\gg(\Sigma_1,z,x)&=(i/2)R_3(z)^{-1/2}\big(\cI_m z+(1/2)\cQ_1(x)
+c_1\cI_m\big), \quad z\in\bbC\backslash\Sigma_1, \; x\in\bbR. \lb{5.36}
\end{align}
with
\begin{equation}
c_1 =-(E_0+E_1+E_2)/2. \lb{5.37}
\end{equation}
By \eqref{5.36}, \eqref{2.54}, \eqref{2.55}, and \eqref{2.65},
$\cQ_1$ is locally absolutely continuous on $\bbR$. By
\eqref{2.66} and \eqref{2.67}, $\cQ'_1$ is locally
absolutely continuous on $\bbR$. Iterating this procedure, using 
\eqref{2.65}--\eqref{2.67b}, one infers inductively that 
\begin{equation}
\cQ_1\in C^\infty(\bbR). \lb{5.38}
\end{equation} 
\eqref{2.69} and \eqref{5.36} then yield
\begin{align}
\gh(\Sigma_1,z,x)&=(i/2)R_3(z)^{-1/2}\big(\cI_m
z^2+(-(1/2)\cQ_1(x) +c_1\cI_m)z \no \\
& \quad +(1/4)\cQ''_1(x)-(1/2)\cQ_1(x)^2
-c_1\cQ_1(x)\big), \lb{5.39} \\
& \hspace*{4.1cm}  z\in\bbC\backslash\Sigma_1, \; x\in\bbR, \no
\end{align}
and \eqref{5.36}, \eqref{2.66a},  and \eqref{2.67a} prove
\begin{equation}
\gg_{p,0}(\Sigma_1,z,x)=-(i/8)R_3(z)^{-1/2}\big(\cQ'_1(x)
+\cC_p\big), \quad p=1,2, \lb{5.40} 
\end{equation}
for some constant matrices $\cC_p\in\bbC^{m\times m}$, $p=1,2$. Insertion
of \eqref{5.36} and \eqref{5.40} into \eqref{2.63}, \eqref{2.64}, taking
into account the asymptotics \eqref{2.50}, \eqref{2.51} of the
Weyl--Titchmarsh matrices $\cM_\pm(\Sigma_1,z,x)$ associated with
$\cQ_1$ then shows $\cC_p=0$, $p=1,2$, and hence
\begin{equation}
\gg_{p,0}(\Sigma_1,z,x)=-(i/8)R_3(z)^{-1/2}\cQ_1'(x), \quad
p=1,2, \; z\in\bbC\backslash\Sigma_1, \; x\in\bbR. \lb{5.41} 
\end{equation}
By \eqref{2.59}, \eqref{5.36}, and \eqref{5.41}, $\cQ_1(x)$ and 
$\cQ'_1(x)$ commute and hence one obtains for each
$x\in\bbR$,
\begin{equation}
[\cQ_1^{(r)}(x),\cQ_1^{(s)}(x)]=0 \text{ for all
$r,s\in\bbN_0$}. \lb{5.42}
\end{equation}
Next, multiplying \eqref{2.62} by $R_3(z)$ and collecting the coefficients
of $z^k$, $0\leq k\leq 2$, yields
\begin{align}
&(1/4)\cQ''_1(x)-(3/4)\cQ_1(x)^2-c_1\cQ_1(x)
+d_1\cI_m=0, \lb{5.43} \\
&((1/4)\cQ''_1(x)-(1/2)\cQ_1(x)^2-c_1\cQ_1(x)) ((1/2)\cQ_1(x)
+c_1\cI_m) \no \\
&-(1/16)\cQ'_1(x)^2 +E_0E_1E_2\cI_m=0, \lb{5.44}
\end{align}
with 
\begin{equation}
d_1=c_1^2-\sum_{\substack{k_1,k_2=0\\k_1<k_2}}^2 E_{k_1}E_{k_2}. \lb{5.45}
\end{equation}
Eliminating $\cQ_1''(x)$ in \eqref{5.43}, \eqref{5.44} finally
yields
\begin{equation}
\cQ'_1(x)^2=-16R_3(-(1/2)\cQ_1(x)-c_1\cI_m). \lb{5.46}
\end{equation} 
Since $\cQ_1(x)$ is a self-adjoint $m\times m$ matrix, we may write
\begin{equation}
\cQ_1(x)=\sum_{k=1}^N q_k(x)\cP_k(x), \lb{5.47}
\end{equation}
where $q_k(x)$ and $\cP_k(x)$ denote the eigenvalues and
corresponding self-adjoint spectral projections of $\cQ_1(x)$, that is, 
\begin{equation}
\cP_k(x)\cP_\ell(x)=\delta_{k,\ell}\cP_\ell(x),
\quad \sum_{k=1}^N \cP_k(x)=\cI_m. \lb{5.48}
\end{equation}
Introducing $\cF_1(\Sigma_1,z,x)$ by 
\begin{equation}
\cF_1(\Sigma_1,z,x)=z\cI_m +(1/2)\cQ_1(x)+c_1\cI_m, \lb{5.49}
\end{equation}
this implies 
\begin{equation}
\cF_1(\Sigma_1,z,x)=\sum_{k=1}^N (z-\mu_k(x))\cP_k(x),
\quad \mu_k(x)=-(1/2)q_k(x)-c_1. \lb{5.50}
\end{equation}
Since by \eqref{4.4},
\begin{equation}
R_3(\lambda)^{1/2}=|R_3(\lambda)|^{1/2}\begin{cases} -1, & \lambda\in
(E_0,E_1), \\ 1, & \lambda\in (E_2,\infty), \end{cases}
\end{equation}
and
\begin{equation}
R_3(\lambda)^{-1/2}\cF_1(\Sigma_1,\lambda,x)>0, \quad
\lambda\in\Sigma_1^o,
\end{equation}
one concludes that for fixed $x\in\bbR$ and all $g\in\bbC^m$, 
$(g,\cF_1(\Sigma_1,\lambda,x)g)_{\bbC^m}$ changes sign for
$\lambda\in [E_1,E_2]$. Thus,
\begin{equation}
\mu_k(x)\in [E_1,E_2], \quad 1\leq k \leq N, \lb{5.50a}
\end{equation}
in accordance with \eqref{4.26}. A comparison with \eqref{4.22} then
yields
\begin{align}
\Gamma_k(\Sigma_1,x)&=-i\lim_{z\to\mu_k(x)}
(z-\mu_k(x))R_3(z)^{1/2}
\cF_1(\Sigma_1,z,x)^{-1} \no \\
&=-iR_3(\mu_k(x))^{1/2}\cP_k(x). \lb{5.51}
\end{align}
Hence, given a sequence of signs, 
\begin{equation}
\varepsilon_k(x)\in\{1,-1\}, \quad 1\leq k\leq N, \lb{5.51b}
\end{equation}
and temporarily assuming 
\begin{equation}
\mu_k(x)\in (E_1,E_2), \quad 1\leq k \leq N, \lb{5.51a}
\end{equation}
one computes
\begin{align}
\cG_{1,0}(\Sigma_1,z,x)&=\bigg(\sum_{k=1}^N
\f{\varepsilon_k(x)}{z-\mu_k(x)}
\Gamma_k(\Sigma_1,x)\bigg)\cF_1(\Sigma_1,z,x) \no \\
&=-i\sum_{k=1}^n \varepsilon_k(x)R_3(\mu_k(x))^{1/2}
\Gamma_k(\Sigma_1,x) \no \\
&=-i R_3(-(1/2)\cQ_1(x)-c_1\cI_m)^{1/2}. \lb{5.52}
\end{align}
Here the choice of the matrix square root in \eqref{5.52} is a
direct consequence of the choice of signs $\varepsilon_k (x)$,
$1\leq k\leq N$. By equation \eqref{5.46}, one obtains 
\begin{equation}
\cG_{1,0}(\Sigma_1,z,x)=-(1/4)\cQ'_1(x) \lb{5.53}
\end{equation}
in accordance with \eqref{5.40}. In particular, $\cF_1(\Sigma_1,z,x)$ in
\eqref{5.49} and $\cG_{1,0}(\Sigma_1,\lambda, x)$ in \eqref{5.53} are of
the form \eqref{5.7} and \eqref{5.8}. Finally, the temporary restriction
\eqref{5.51a} can be removed by continuity. Summing up, 
\begin{equation}
\cQ_1(x)=\cQ_{\Sigma_1}(x), \quad x\in\bbR, \lb{5.54}
\end{equation}
with $\cQ_{\Sigma_1}$ constructed as in Section \ref{s3} given
some $\cF_{1,\sigma_1}(z,x_0)$ and some choice of signs 
$\varepsilon_k(x_0)\in\{1,-1\}$, $1\leq k\leq N$ (cf.\ also 
\eqref{5.7}--\eqref{5.12}). 
\end{proof} 

We note that the case $\ell=0$ in Theorem \ref{t4.2} was originally 
treated in \cite{CGHL00} using a (matrix-valued) trace formula 
approach. 

Combining Theorems \ref{t4.1} and \ref{t4.2} proves Theorem \ref{t1.5}. 

Since (self-adjoint) periodic potentials $\cQ$ which lead to
Schr\"odinger operators with uniform maximum spectral multiplicity are
reflectionless in the sense of Definition \ref{d2.7} as shown in
\cite{CGHL00}, Theorem \ref{t1.5a} yields the proper matrix
generalizations of Borg's and Hochstadt's results, Theorem \ref{t1.1} 
and \ref{t1.3}.

\bigskip
\noindent {\bf Acknowledgements.} 
It is a great pleasure to dedicate this paper to Jerry Goldstein 
and Rainer Nagel on the occasion of their 60th birthdays. In particular,
F.\ G. gratefully acknowledges Jerry's friendship and constant support
throughout the years. 
 
We thank Robert Carlson, Mark Malamud, Fedor Rofe-Beketov, and  Barry
Simon for many helpful discussions on the material presented in this
paper.

Research leading to this paper started in 1999 when E.\ B.\ and F.\ G.\
were supported in part by the CRDF grant UM1-325. 
  


\end{document}